\documentclass[11pt]{amsart}

\usepackage{amsmath,amssymb,amsthm,mathtools}
\usepackage{booktabs}
\usepackage[margin=1.2in]{geometry}
\usepackage{hyperref}

\makeatletter
\renewcommand{\subsubsection}{\@startsection{subsubsection}{3}%
  {\z@}{.8\baselineskip\@plus.2\baselineskip}%
  {.25\baselineskip}{\normalfont\bfseries}}
\makeatother

\newenvironment{tablelayout}
  {%
   \setlength{\abovedisplayskip}{12pt plus 3pt minus 2pt}%
   \setlength{\belowdisplayskip}{12pt plus 3pt minus 2pt}%
   \setlength{\abovedisplayshortskip}{\abovedisplayskip}%
   \setlength{\belowdisplayshortskip}{\belowdisplayskip}}
  {}

\numberwithin{equation}{section}

\newtheorem{theorem}{Theorem}[section]
\newtheorem{lemma}[theorem]{Lemma}
\newtheorem{proposition}[theorem]{Proposition}
\newtheorem{corollary}[theorem]{Corollary}
\theoremstyle{definition}
\newtheorem{definition}[theorem]{Definition}
\theoremstyle{plain}
\newtheorem{maintheorem}{Theorem}

\newtheorem{mainproposition}[maintheorem]{Proposition}

\newcommand{\R}{\mathbb R}
\newcommand{\C}{\mathbb C}
\newcommand{\KW}{\mathcal K_W}
\newcommand{\GW}{\mathcal G_W}
\newcommand{\CW}{\mathcal C_W}
\newcommand{\Eh}{\widehat E}
\newcommand{\Vh}{\widehat V}
\newcommand{\Cl}{\mathrm{Cl}}
\DeclareMathOperator{\tr}{tr}
\DeclareMathOperator{\Scal}{Scal}
\DeclareMathOperator{\Hol}{Hol}
\DeclareMathOperator{\End}{End}
\DeclareMathOperator{\ev}{ev}

\title{Generalised Killing spinors and Euclidean hypersurfaces}
\author{Samuel Forster}
\author{Marie-Am\'elie Lawn}
\address{Department of Mathematics, Imperial College London}
\date{\today}
\subjclass[2020]{53C27, 53C42, 53C29}
\keywords{Generalised Killing spinors, hypersurfaces,
spin geometry, holonomy}

\begin{document}

\begin{abstract}

We study when the number of spinors satisfying the same
generalised Killing equation for a symmetric endomorphism
field $W$ forces local realisation as a Euclidean
hypersurface. Let $N=2^{\lfloor n/2\rfloor}$ be the rank
of the complex spinor bundle. We show that in dimension
$n\geq4$, if the complex dimension of the solution space
exceeds $N/2$, then the metric and $W$ satisfy the
Gauss--Codazzi equations. This is proved in the smooth category by combining
the curvature defects into a single bivector and
bounding the dimension of its kernel on spinors,
using algebraic curvature symmetries to handle
exceptional cases. Examples
obtained from non-flat four-dimensional hyperk\"ahler
manifolds show that this threshold is sharp for every
$n\geq4$. In dimension three, the Heisenberg group carries
the maximal number of common solutions while failing
the Euclidean Gauss equation. Further, in the real analytic
category, we use existing extension theorems for Einstein
metrics and parallel spinors to extend all common
solutions to parallel spinors on a single local Ricci-flat
ambient manifold. Using Wang's counts of parallel spinors, we obtain
an independent proof of the sharp threshold for
real analytic data.

\end{abstract}

\maketitle

\section{Introduction}

The classical hypersurface problem asks when a Riemannian
manifold $(M^n,g)$, equipped with a symmetric endomorphism
field $W\in\Gamma(\End(TM))$, admits a local isometric
immersion into $\R^{n+1}$ with shape operator $W$.
The fundamental theorem of hypersurfaces gives a complete
local answer to this problem: such an immersion exists
if and only if $g$ and $W$ satisfy the Gauss--Codazzi
equations. Explicitly, these can be written as follows.
\begin{align}
R(X,Y)U
&= g(W(Y),U)W(X)-g(W(X),U)W(Y),
\label{eq:gauss}\\
(\nabla_X W)(Y)
&= (\nabla_Y W)(X).
\label{eq:codazzi}
\end{align}

Spinors provide another way of characterising such
hypersurface data. This approach was first developed by
Sullivan \cite{Sullivan1986} for minimal surfaces in $\R^3$,
which he described and reconstructed using spinors,
and was subsequently extended by Kusner--Schmitt
\cite{KusnerSchmitt1996} to general conformal immersions.
Friedrich \cite{Friedrich1998}
gave an intrinsic characterisation. He showed that
restricting a non-zero parallel spinor of the ambient
space to an oriented immersed surface
$M^{2} \rightarrow \R^3$, with a suitable identification of
the spinor bundles, gives a non-zero solution of constant
length to the Dirac equation with potential
\[
D\psi = H\psi,
\]
where $H$ is the mean curvature function. Conversely, a non-zero solution of constant length
determines a smooth symmetric endomorphism field $W$
through
\[
g(W(X),Y)
=-\frac{2}{|\psi|^2}
\operatorname{Re}
\langle\nabla_X^{\Sigma M}\psi,Y\cdot\psi\rangle,
\]
where $\langle\cdot,\cdot\rangle$ is the Hermitian
spinor inner product. This field satisfies
$\tr(W)=2H$, and $\psi$ satisfies the following
equation, now called the generalised Killing
spinor equation:
\begin{equation}
  \nabla_X^{\Sigma M}\psi
  =-\frac12W(X)\cdot\psi,
  \qquad \forall X\in\Gamma(TM).
  \label{eq:gks}
\end{equation}
In dimension two, the integrability conditions of
this equation imply that $g$ and $W$ satisfy the
Gauss--Codazzi equations. Hence one obtains a local
isometric immersion with shape operator $W$.

Morel \cite{Morel2005} extended the surface correspondence
to $S^3$ and $H^3$. He also showed that, in dimension three,
a single non-zero generalised Killing spinor gives a local
isometric immersion into $\R^4$ with shape operator $W$
if $W$ is assumed to be a Codazzi tensor.
Lawn--Roth \cite{LawnRoth2010} obtained a characterisation
using two spinors, without imposing the Codazzi equation
separately. In the Euclidean case, their criterion requires
non-zero spinors $\psi_+$ and $\psi_-$, one for $W$ and one
for $-W$, satisfying
\[
\nabla_X^{\Sigma M}\psi_\pm
= \mp\frac12 W(X)\cdot\psi_\pm,
\qquad \forall X\in\Gamma(TM).
\]
Their existence forces the Gauss--Codazzi equations and
hence gives a local isometric immersion into $\R^4$.

In arbitrary dimension and codimension,
Bayard--Lawn--Roth \cite{BayardLawnRoth2017} gave a spinorial
representation of submanifolds in Riemannian space forms.
In their Euclidean formulation, the twisted spinor bundle
is constructed using left multiplication on the full real
Clifford algebra. One solution is enough to recover the immersion,
provided its values in a spin frame lie in the ambient
spin group. Here we work with the intrinsic spinor
bundle, and impose no such algebraic condition on
solutions of \eqref{eq:gks}.

A related question asks whether a generalised Killing
spinor extends to a parallel spinor on an ambient manifold.
B\"ar--Gauduchon--Moroianu \cite{BGM2005} answered this
affirmatively when $W$ is a Codazzi tensor, by constructing
an explicit generalised cylinder carrying such an
extension. In the Riemannian setting, the resulting
ambient metric is Ricci-flat.
Ammann--Moroianu--Moroianu \cite{AMM2013} removed the
Codazzi assumption in the real analytic category:
if $g$ and $W$ are real analytic, a single non-zero
generalised Killing spinor extends to a parallel spinor
on a local Ricci-flat ambient manifold realising $(g,W)$
as hypersurface data. They also showed that this extension
need not exist for general smooth data.
In ambient dimensions six to eight, related extension
results are expressed through hypo structures
\cite{ContiSalamon2007} and Hitchin's evolution equations
for differential forms \cite{Hitchin2001}.

For Euclidean hypersurfaces in higher dimensions, the
existence of a single generalised Killing spinor in the
intrinsic spinor bundle does not by itself force the
Gauss--Codazzi equations. We call the differences between
the two sides of \eqref{eq:gauss} and \eqref{eq:codazzi}
the Gauss and Codazzi defects, respectively. In the
curvature identity associated with equation~\eqref{eq:gks},
the Gauss defect appears as a degree-two Clifford term
and the Codazzi defect as a degree-one term, both acting
on the same spinor. The two terms may cancel on a spinor
even when neither is zero. For $n\geq4$, we know that
taking the maximal possible number of independent
solutions forces the Gauss--Codazzi equations
(Lemma~\ref{lem:maximal}). It is natural to ask how many
solutions we actually need. We look for the smallest
number of independent solutions for one fixed $W$ that
forces these equations. Dimension three is different:
$\mathfrak{so}(4)$ is not simple, and its action on either
half-spin module has a non-zero kernel. As explained
after Lemma~\ref{lem:maximal}, this allows non-zero
curvature defects to annihilate the whole intrinsic
spinor fibre. Proposition~\ref{prop:C} shows that even
maximal multiplicity need not force the Gauss--Codazzi
equations.

For a symmetric endomorphism field
$W$, we define the vector space of its associated generalised Killing
spinors by
\[
\mathcal K_W(M):=
\left\{
\psi\in\Gamma(\Sigma M):
\nabla_X^{\Sigma M}\psi
=-\tfrac12 W(X)\cdot\psi
\quad\forall X\in\Gamma(TM)
\right\}.
\]
The complex rank of the spinor bundle is
\[
N:=\operatorname{rank}_\C\Sigma M
=2^{\lfloor n/2\rfloor}.
\]
In what follows, what matters is not the number of
solutions gathered over different fields $W$, but the
number of common solutions for a single fixed field $W$.

\begin{maintheorem}[Sharp threshold in the smooth category]
\label{thm:A}
Let $(M^n,g)$ be a smooth, connected Riemannian spin
manifold of dimension $n\geq4$, and let $W$ be a smooth
field of symmetric endomorphisms of the tangent bundle
$TM$. If
\[
\dim_\C\mathcal K_W(M)>N/2,
\]
then the Euclidean Gauss and Codazzi equations are
satisfied. Consequently, $(g,W)$ is locally realised as
the induced metric and shape operator of an isometric
immersion as a hypersurface in $\mathbb R^{n+1}$.
Equivalently, the hypothesis is
$\dim_\C\KW(M)\geq2^{\lfloor n/2\rfloor-1}+1$.
\end{maintheorem}

The sharpness examples use the restriction of ambient
parallel spinors. The Ricci-flat extensions constructed in \cite{AMM2013}
need not be flat, which suggests looking for non-flat
Ricci-flat manifolds with as many parallel spinors as
possible. Restricting them to suitable hypersurfaces
produces many solutions for one shape operator $W$;
we choose the hypersurfaces so that a non-zero
tangential component of the ambient curvature remains
as a defect in the Euclidean Gauss equation.

\begin{maintheorem}[Sharpness of the threshold]
\label{thm:B}
In every dimension $n\geq4$ there exist smooth data
$(M^n,g,W)$ with
\[
\dim_\C\KW(M) = N/2
\]
which do not satisfy the Euclidean Gauss equation.
These examples can also be taken real analytic.
\end{maintheorem}

For ordinary Killing spinors, the equation is
\[
\nabla_X\psi=\lambda X\cdot\psi,
\]
where $\lambda\neq0$ is a fixed real constant.
Alekseevsky--Cort\'es show that more than $3N/8$
independent solutions for the same $\lambda$
force constant sectional curvature $4\lambda^2$
\cite[Theorem~8]{AlekseevskyCortes2009}.
Thus the metric is locally that of a round sphere
of radius $1/(2|\lambda|)$.

In our notation, this is the special case
$W=-2\lambda\,\mathrm{id}$. Codazzi holds
automatically, and the Gauss equation reduces to
\[
R(X,Y)Z
=4\lambda^2\bigl(g(Y,Z)X-g(X,Z)Y\bigr).
\]
Their result therefore gives a local Euclidean
hypersurface immersion with the prescribed shape
operator $W$, realised by the round sphere with
the appropriate choice of normal.
Theorem~\ref{thm:A} gives the same conclusion
in this case, but requires more than $N/2$
solutions. The larger threshold is necessary
when we allow arbitrary symmetric $W$:
our sharpness examples show that it cannot
be lowered in that generality.

Wang's count \cite{Wang1989} gives a simple way to
see why dimension three is different. A non-flat
hyperk\"ahler four-manifold carries two independent
parallel spinors in the same half-spin bundle.
Restricting these to a suitably oriented
three-dimensional hypersurface gives two common
generalised Killing spinors, the maximal number in
this dimension. We can choose the tangent space at
one point to contain a plane of non-zero ambient
sectional curvature, so the Euclidean Gauss equation
fails. We give an explicit example on the Heisenberg
group below.

\begin{mainproposition}[Failure in dimension three]
\label{prop:C}
The three-dimensional Heisenberg group carries a
left-invariant metric $g$ and a left-invariant symmetric
endomorphism field $W$ with
\[
\dim_\C\KW(M) = 2 = N
\]
which do not satisfy the Euclidean Gauss equation.
Its generalised Killing spinors are among Artacho's
examples of left-invariant spinors \cite{Artacho2025}.
\end{mainproposition}

Another consequence of Theorem~\ref{thm:A} is a local
multiplicity gap. Let $(M^n,g,W)$ be as in
Theorem~\ref{thm:A}, and let $U\subset M$ be a connected,
simply connected open set. Then
\begin{equation}
\dim_\C\KW(U) \in \{0,1,\dots,N/2\}\cup\{N\}.
\label{eq:gap}
\end{equation}

In the smooth proof, we adjoin a formal unit normal to
the tangent space and combine the two defects into a
single real bivector. Estimating the kernel of its action
on spinors gives the conclusion immediately in most
dimensions. In the remaining dimensions $7,11,15$, we use
the fact that all curvature values annihilate the same
subspace of spinors. This common annihilator can only
be zero, the line spanned by a K\"ahler form, or the
three-dimensional space associated with a quaternionic
structure. The pair symmetry of the Gauss defect and the
first Bianchi identity force the curvature values lying
in the latter two to vanish as well.

For real analytic data, we apply Koiso's metric extension
theorem \cite{Koiso1981} and the parallel spinor extension
theorem of Ammann--Moroianu--Moroianu \cite{AMM2013} to the
whole solution space for fixed data $(g,W)$. All common
solutions extend to a single ambient metric, on which
we can then count parallel spinors. We therefore obtain the sharp multiplicity threshold by two independent arguments: a pointwise algebraic
proof for smooth data, and a proof for real analytic
data using Ricci-flat extensions and holonomy.

For a Riemannian spin manifold $Z$, write $\mathcal P(Z)$
for the complex vector space of parallel spinors and,
when $\dim Z$ is even, $\mathcal P^\pm(Z)$ for the
subspaces of positive and negative chirality. For an
oriented hypersurface $M^n\subset Z^{n+1}$ with $n$ odd,
$\tau\in\{+,-\}$ denotes the half-spin choice compatible
with our intrinsic Clifford multiplication and the
hypersurface identification in \eqref{eq:spingauss}.

\begin{maintheorem}[Analytic extension and parallel spinor counts]
\label{thm:D}
Let $(M^n,g)$ be a connected real analytic Riemannian
spin manifold, and let $W$ be a real analytic symmetric
endomorphism field with $\KW(M)\neq0$.
On a sufficiently small connected, simply connected
neighbourhood $U$ of each point, the following hold.
\begin{enumerate}
\item[(i)]
For some $\varepsilon>0$, there is a real analytic
Ricci-flat metric $\bar g=dt^2+g_t$ on
$Z=(-\varepsilon,\varepsilon)\times U$
whose initial slice has induced metric $g$ and shape
operator $W$ with respect to $\partial_t$.
This metric is uniquely determined near the initial
slice by $(g,W)$ in these Gaussian coordinates.

\item[(ii)]
Every spinor in $\KW(U)$ extends uniquely to a parallel
spinor for this same ambient metric, in the appropriate
spinor bundle. Restriction induces isomorphisms
\[
\begin{aligned}
\mathcal P(Z)&\cong\KW(U),
&&n\text{ even},\\
\mathcal P^\tau(Z)&\cong\KW(U),
&&n\text{ odd}.
\end{aligned}
\]

\item[(iii)]
If $n\geq4$ and $\bar g$ is non-flat, then
\[
\dim_\C\KW(U)\leq N/2.
\]
For every $n\geq4$, there are examples with $\bar g$
non-flat and $\dim_\C\KW(U)=N/2$.
\end{enumerate}
\end{maintheorem}

For the non-flat upper bound, we use the local
de Rham decomposition and Wang's parallel-spinor
counts \cite{Wang1989,MoroianuSemmelmann2000}.
Alekseevsky--Cort\'es combine these tools with
product formulas to obtain a full-spin flatness
criterion
\cite[proof of Theorem~4(ii)]{AlekseevskyCortes2009}.
We apply this approach to the common Ricci-flat
extension. In odd hypersurface dimensions, restriction
uses only one ambient half-spin bundle, so we also
need to count the parallel spinors of each chirality
separately.

Section~\ref{sec:conventions} fixes conventions and the
curvature defects, and Section~\ref{sec:examples} gives
the examples. Section~\ref{sec:smooth} proves the smooth
threshold, and Section~\ref{sec:analytic} treats the
analytic extension and the count estimates.
Appendix~\ref{app:LO} proves the combinatorial bound on
bivector kernels used in Section~\ref{sec:smooth}, together
with the characterisation of its exceptional cases.

\section{Intrinsic generalised Killing data and
curvature defects}
\label{sec:conventions}

\subsection{The common equation and the modified connection}

We take the curvature and Clifford multiplication
conventions
\begin{equation}
R(X,Y) = [\nabla_X,\nabla_Y]-\nabla_{[X,Y]},
\qquad
X\cdot Y+Y\cdot X = -2g(X,Y).
\label{eq:conventions}
\end{equation}
The spinor bundle carries the standard Hermitian inner
product, and Clifford multiplication by real tangent
vectors is skew-Hermitian. For an isometric hypersurface
immersion $F:M\to Z$, the shape operator is defined by
\[
dF(W(X)) = -\bar\nabla_{dF(X)}\nu,
\]
where $\bar\nabla$ is the ambient Levi-Civita connection
and $\nu$ is the chosen unit normal. Identifying the
tangent space with its image,
$W(X) = -\bar\nabla_X\nu$.

\begin{definition}
The modified connection associated with $W$ is defined by
\begin{equation}
\nabla^W_X := \nabla^{\Sigma M}_X + \frac12 W(X)\cdot.
\label{eq:modconn}
\end{equation}
\end{definition}

The correction term is skew-Hermitian, so this connection
preserves the Hermitian inner product, and $\KW(M)$ is
exactly the space of $\nabla^W$-parallel sections. Hence
the lengths of solutions and their inner products
are constant. On connected $M$, a $\nabla^W$-parallel section that
vanishes at one point vanishes identically, since
parallel transport determines its value along every
curve from that point.
Hence at each $p\in M$ the evaluation map
\begin{equation}
\ev_p:\KW(M)\longrightarrow\Sigma_pM,
\qquad \psi\longmapsto\psi(p)
\label{eq:eval}
\end{equation}
is injective, and
$\dim_\C\KW(M)\leq N=\operatorname{rank}_\C\Sigma M$.
In particular, common solutions that are globally
linearly independent are linearly independent at every
point.

\subsection{The Gauss and Codazzi defects}

We define the following tensors.
\begin{align}
(W\wedge W)(X,Y)U
&:= g(W(Y),U)W(X)-g(W(X),U)W(Y),
\label{eq:wedge}\\
\GW(X,Y)
&:= R(X,Y)-(W\wedge W)(X,Y),
\label{eq:gaussdefect}\\
\CW(X,Y)
&:= (\nabla_XW)(Y)-(\nabla_YW)(X),
\label{eq:codazzidefect}\\
\GW(X,Y,U,V)
&:= g(\GW(X,Y)U,V).
\label{eq:gauss4}
\end{align}
For fixed $X,Y$, $\GW(X,Y)$ is a skew endomorphism,
whereas $\CW(X,Y)$ is a tangent vector.
By the fundamental theorem of hypersurfaces,
$\GW=\CW=0$ is necessary and sufficient for $(g,W)$
to be locally realised as Euclidean hypersurface data.
For a proof using generalised cylinders, see
\cite{BGM2005}.

\begin{proposition}
The curvature of the modified connection is
\begin{equation}
R^{\nabla^W}_{X,Y}
= R^{\Sigma M}(X,Y)+\frac12\CW(X,Y)\cdot
+\frac14\bigl[W(X)\cdot,W(Y)\cdot\bigr].
\label{eq:modcurv}
\end{equation}
\end{proposition}

\begin{proof}
By compatibility with Clifford multiplication, for
every spinor section $\psi$,
\[
\nabla^{\Sigma M}_X(W(Y)\cdot\psi)
= \nabla_X(W(Y))\cdot\psi
+W(Y)\cdot\nabla^{\Sigma M}_X\psi.
\]
Expand $R^{\nabla^W}_{X,Y}\psi$ using
\eqref{eq:modconn}. The mixed terms containing first
derivatives of $\psi$ cancel when we subtract the
expression with $X$ and $Y$ exchanged. The derivatives
of $W$, together with the correction from
$\nabla^W_{[X,Y]}$, give
\[
\frac12\bigl(
\nabla_X(W(Y))-\nabla_Y(W(X))-W([X,Y])
\bigr)
= \frac12\CW(X,Y),
\]
using the torsion-free identity
$[X,Y]=\nabla_XY-\nabla_YX$.
The terms quadratic in $W$ give the commutator
in \eqref{eq:modcurv}.
\end{proof}

With respect to an orthonormal frame $\{e_i\}$, the spin
curvature is
\[
R^{\Sigma M}(X,Y)\psi
= \frac14\sum_{i,j}
g(R(X,Y)e_i,e_j)\,e_i\cdot e_j\cdot\psi.
\]
On the other hand, expanding
$W(X)=\sum_i g(W(X),e_i)e_i$ and similarly for $W(Y)$
gives
\[
\frac14\sum_{i,j}
g((W\wedge W)(X,Y)e_i,e_j)\,
e_i\cdot e_j\cdot\psi
= \frac14\bigl(
W(Y)\cdot W(X)-W(X)\cdot W(Y)
\bigr)\cdot\psi.
\]
For any spinor field $\psi$, the preceding identities
and antisymmetry in $i,j$ give
\begin{equation}
R^{\nabla^W}_{X,Y}\psi
= \frac12\sum_{i<j}
\GW(X,Y,e_i,e_j)\,e_i\cdot e_j\cdot\psi
+\frac12\CW(X,Y)\cdot\psi.
\label{eq:intid}
\end{equation}
If $\psi\in\KW(M)$, then $\psi$ is parallel for
$\nabla^W$, so $R^{\nabla^W}_{X,Y}\psi=0$.

\subsection{Adjoining a formal normal}
\label{subsec:formalnormal}

We now turn the degree-one Codazzi term into a bivector
by adjoining one extra direction. Fix a point $p$
and set $E=T_pM$ and $\Eh=E\oplus\R e_0$, where $e_0$
is a formal unit vector orthogonal to $E$. The entire
curvature will then act through $\Lambda^2\Eh$.
We identify bivectors with degree-two Clifford elements
by $u\wedge v=\frac12(u\cdot v-v\cdot u)$. The map
$X\mapsto X\cdot e_0$ satisfies
\[
(X\cdot e_0)(Y\cdot e_0)
+(Y\cdot e_0)(X\cdot e_0)
=-2\langle X,Y\rangle,
\]
and so extends to an isomorphism of complex Clifford
algebras
\begin{equation}
\alpha:\Cl_\C(E)\longrightarrow\Cl^0_\C(\Eh),
\qquad \alpha(X)=X\cdot e_0.
\label{eq:alpha}
\end{equation}
Indeed, the image generates the even subalgebra, and
both algebras have complex dimension $2^n$. Moreover, $(X\cdot e_0)(Y\cdot e_0)=X\cdot Y$, so
$\alpha$ is the identity on $\Lambda^2E$, where we
realise $\Lambda^2E$ as a subspace of both Clifford
algebras by the above identification.

After choosing a spin frame at $p$, identify $\Sigma_pM$
with $\Delta_n$. For even $n$, $\Cl_\C(E)$ has a unique
irreducible module up to isomorphism, realised on
$\Delta_{n+1}$ through $\alpha$. For odd $n$, the two
half-spin modules realise its two inequivalent
irreducible representations; let $\tau$ denote the
unique choice isomorphic to the chosen intrinsic
Clifford module. There is therefore an isomorphism
\begin{equation}
\begin{aligned}
\iota:\Delta_n&\longrightarrow\Delta_{n+1},
&&n\text{ even},\\
\iota:\Delta_n&\longrightarrow\Delta^\tau_{n+1},
&&n\text{ odd},
\end{aligned}
\label{eq:iota}
\end{equation}
satisfying
\[
\iota(a\cdot\xi)=\alpha(a)\cdot\iota(\xi).
\]
In either case the target has dimension $N$. Set
\[
\Vh_p:=\{\iota(\psi(p)):\psi\in\KW(M)\}.
\]
By \eqref{eq:eval},
$\dim_\C\Vh_p=\dim_\C\KW(M)$. Under $\alpha$, the
degree-two term in \eqref{eq:intid} is unchanged,
while $\CW(X,Y)$ is sent to
$\CW(X,Y)\cdot e_0=\CW(X,Y)\wedge e_0$. Define the
real bivector
\begin{equation}
\Omega(X,Y)
:=\sum_{i<j}\GW(X,Y,e_i,e_j)\,e_i\wedge e_j
+\CW(X,Y)\wedge e_0.
\label{eq:Omega}
\end{equation}
Then, for every $\xi\in\Sigma_pM$,
\[
\iota\bigl(R^{\nabla^W}_{X,Y}\xi\bigr)
=\frac12\Omega(X,Y)\cdot\iota(\xi).
\]
In particular,
\begin{equation}
\Omega(X,Y)\cdot\Vh_p=0.
\label{eq:annihil}
\end{equation}
The induced inner product on bivectors is characterised by
\[
\langle u\wedge v,U\wedge V\rangle
=\langle u,U\rangle\langle v,V\rangle
-\langle u,V\rangle\langle v,U\rangle.
\]
Thus, for $U,V\in E$,
\begin{equation}
\begin{aligned}
\langle\Omega(X,Y),U\wedge V\rangle
&=\GW(X,Y,U,V),\\
\langle\Omega(X,Y),U\wedge e_0\rangle
&=g(\CW(X,Y),U).
\end{aligned}
\label{eq:Omega-components}
\end{equation}
These component identities show that
$\Omega(X,Y)=0$ for all $X,Y\in E$ if and only if
$\GW=\CW=0$ at $p$. This construction is an algebraic
operation at a single point, and does not assume the
existence of an ambient manifold or of an actual
immersion.

We have reduced the Gauss--Codazzi equations to showing
that $\Omega(X,Y)=0$ for all $X,Y\in E$. Each of these
bivectors annihilates the values of all common solutions.
We first consider the maximal possible multiplicity,
when these values span the whole spinor fibre, and show
that this is sufficient to force both defects to vanish.

\begin{lemma}[Maximal multiplicity]
\label{lem:maximal}
Let $(M^n,g)$ be a connected smooth Riemannian spin
manifold of dimension $n\geq4$, and let $W$ be a smooth
symmetric endomorphism field. If
\[
\dim_\C\KW(M)=N,
\]
then $g$ and $W$ satisfy the Euclidean Gauss--Codazzi
equations.
\end{lemma}

\begin{proof}
Fix $p\in M$. Injectivity of evaluation implies that
the values of the common solutions span $\Sigma_pM$.
Using \eqref{eq:iota}, let
\[
\mathrm{cl}:\Lambda^2\Eh
\longrightarrow\End_\C(\Sigma_pM)
\]
denote Clifford action transported to the intrinsic
spinor fibre. By \eqref{eq:annihil}, every curvature
bivector annihilates all solution values, so
\[
\mathrm{cl}(\Omega(X,Y))=0
\qquad\text{for all }X,Y\in E.
\]

We need to check that only the zero bivector acts
as zero on the whole module; that is, that
$\mathrm{cl}$ is faithful. Inside $\Cl^0(\Eh)$,
the real bivectors form the Lie algebra
$\mathfrak{spin}(n+1)$ under the Clifford commutator
$[A,B]_{\Cl}:=AB-BA$. Its identification with
$\mathfrak{so}(\Eh)$ is
\begin{equation}
\Lambda^2\Eh\longrightarrow\mathfrak{so}(\Eh),
\qquad
A\longmapsto\bigl(z\longmapsto[A,z]_{\Cl}\bigr).
\label{eq:spin-adjoint}
\end{equation}
Indeed, the Clifford relations give
\[
[u\wedge v,z]_{\Cl}
=2\bigl(\langle u,z\rangle v
        -\langle v,z\rangle u\bigr),
\]
so this map identifies a basis of bivectors with
a basis of skew endomorphisms. It preserves brackets
by the Jacobi identity.

Clifford multiplication gives a representation of
this Lie algebra:
\[
\mathrm{cl}([A,B]_{\Cl})
=\mathrm{cl}(A)\mathrm{cl}(B)
-\mathrm{cl}(B)\mathrm{cl}(A).
\]
Its kernel is therefore an ideal: if
$\mathrm{cl}(A)=0$, then
$\mathrm{cl}([B,A]_{\Cl})=0$ for every bivector $B$.
For $n\geq4$, the Lie algebra
$\mathfrak{spin}(n+1)\cong\mathfrak{so}(n+1)$ is simple,
so this kernel is either zero or the whole algebra.
The latter is impossible, since for orthonormal
$e_1,e_2\in E$,
\[
\mathrm{cl}(e_1\wedge e_2)^2=-\mathrm{id}.
\]
Consequently, $\ker\mathrm{cl}=0$, and hence
$\Omega(X,Y)=0$ for all $X,Y\in E$. Applying
\eqref{eq:Omega-components} at each point gives
$\GW=\CW=0$ throughout $M$.
\end{proof}

For $n=3$, the Lie algebra instead splits as
\[
\mathfrak{spin}(4)
\cong\mathfrak{su}(2)\oplus\mathfrak{su}(2).
\]
Each half-spin representation is faithful on one
summand and zero on the other. Since the intrinsic
spin module is identified with one of these half-spin
modules, a non-zero bivector can act as zero on the
whole intrinsic spinor fibre. Thus maximal multiplicity
does not force $\Omega=0$. The Heisenberg example in
Section~\ref{sec:examples} has maximal multiplicity
but fails the Gauss--Codazzi equations.

\section{Sharpness and the failure in dimension three}
\label{sec:examples}

We obtain generalised Killing spinors by restricting
ambient parallel spinors to a hypersurface. Let
$M^n\subset Z^{n+1}$ be an oriented hypersurface of a
Riemannian spin manifold, with unit normal $\nu$ and
the induced spin structure. At each point of $M$,
the normal $\nu$ replaces the formal vector $e_0$.
The module identification \eqref{eq:iota} then gives
a bundle isomorphism
\[
\iota:\Sigma M\longrightarrow
\begin{cases}
(\Sigma Z)|_M, & n\text{ even},\\
(\Sigma^\tau Z)|_M, & n\text{ odd},
\end{cases}
\]
with the same half-spin choice $\tau$. It satisfies
\[
\iota(X\cdot_M\psi)
=X\cdot_Z\nu\cdot_Z\iota(\psi).
\]
For an ambient spinor $\Psi$ in the indicated bundle,
we define its intrinsic restriction by
\[
\psi:=\iota^{-1}(\Psi|_M).
\]
Under the same identification, the spinorial Gauss
formula \cite{BGM2005} reads
\begin{equation}
\bigl(\bar\nabla^{\Sigma Z}_X\Psi\bigr)\big|_M
=\nabla^{\Sigma M}_X\psi
+\frac12W(X)\cdot_M\psi.
\label{eq:spingauss}
\end{equation}
If $\Psi$ is parallel, the left-hand side vanishes,
so $\psi$ satisfies \eqref{eq:gks}.

For odd $n$, Clifford multiplication by $\nu$ gives
the identification with the opposite half-spin bundle:
\begin{equation}
\begin{aligned}
\iota^{\mathrm{opp}}(\psi)
&:=\nu\cdot_Z\iota(\psi)
\in(\Sigma^{-\tau}Z)|_M,\\
\iota^{\mathrm{opp}}(X\cdot_M\psi)
&=-X\cdot_Z\nu\cdot_Z\iota^{\mathrm{opp}}(\psi).
\end{aligned}
\label{eq:iota-opp}
\end{equation}
Under this identification, the spinorial Gauss
formula has the opposite sign in its $W$ term, since
$\bar\nabla_X\nu=-W(X)$ and
$\nu\cdot_Z W(X)\cdot_Z\nu=W(X)$.
Thus an ambient parallel spinor in $\Sigma^{-\tau}Z$
restricts to a solution of
\[
\nabla_X^{\Sigma M}\psi
=+\frac12W(X)\cdot_M\psi,
\]
the generalised Killing equation for $-W$.
When counting solutions for the prescribed field $W$,
we therefore use the parallel spinors in
$\Sigma^\tau Z$.

Restriction is injective on the space of ambient
parallel spinors. Indeed, a parallel spinor whose
restriction is zero vanishes at a point of $M$,
and hence on the whole connected ambient neighbourhood.

\subsection{Half-rank examples in dimension four and higher}
\label{subsec:halfrank}

Let $X^4$ be a non-flat hyperk\"ahler four-manifold with
holonomy $\mathrm{Sp}(1)$. For a compact example, fix a
K\"ahler class on a K3 surface and take the unique
Ricci-flat K\"ahler metric in that class, given by Yau's
theorem \cite{Yau1978}. With a suitable orientation,
$X$ has two complex linearly independent parallel spinors,
both in the same half-spinor bundle \cite{Wang1989}.
That is,
\[
\dim_\C\mathcal P^+(X)=2,\qquad
\dim_\C\mathcal P^-(X)=0.
\]
For $n\geq5$ consider the product hypersurface
\begin{equation}
M^n = X^4\times S^{n-4}(r)
\subset X^4\times\R^{n-3} = Z^{n+1}.
\label{eq:product}
\end{equation}
The spin structure is induced from the ambient product.
Choosing the outward unit normal in the sphere direction,
\begin{equation}
W|_{TX}=0,\qquad
W|_{TS^{n-4}(r)}=-r^{-1}\,\mathrm{id}.
\label{eq:productW}
\end{equation}
Indeed, the normal is constant in the $X$-directions,
and $\bar\nabla_Y\nu=r^{-1}Y$ in the sphere directions.

Since $\mathcal P^-(X)=0$, the product rules reduce to
\[
\begin{aligned}
\mathcal P(Z)
&\cong\mathcal P^+(X)\otimes\mathcal P(\R^{n-3}),
&&n\text{ even},\\
\mathcal P^\pm(Z)
&\cong\mathcal P^+(X)\otimes\mathcal P^\pm(\R^{n-3}),
&&n\text{ odd}.
\end{aligned}
\]
When $n$ is even, the flat factor has odd dimension
$n-3$, so the number of ambient parallel spinors is
\[
2\cdot2^{(n-4)/2}=2^{n/2-1}=N/2.
\]
When $n$ is odd, $n-3\geq2$ is even and each half-spin
module of the flat factor has dimension $2^{(n-5)/2}$.
Hence the number of parallel spinors in each ambient
half-spinor bundle is
\[
2\cdot2^{(n-5)/2}=2^{(n-3)/2}=N/2.
\]
In either case, a suitable restriction gives
$\dim_\C\KW(M)\geq N/2$.

On the other hand, for $u,v,w\in TX$,
\eqref{eq:productW} gives
\[
\GW(u,v)w = R^X(u,v)w.
\]
Since $X$ is non-flat, the Gauss defect does not vanish
identically. In the four-dimensional case it suffices
to take
\begin{equation}
M^4=X^4\times\{0\}\subset X^4\times\R,\qquad W=0.
\label{eq:slice}
\end{equation}
The two parallel spinors give common solutions, and
again $\GW=R^X\neq0$. To see that sharpness persists
when the generalised Killing spinors are not parallel
(equivalently, when $W$ is not identically zero),
choose $p\in X$ with $R^X_p\neq0$ and replace the slice
by the graph of a function $f$ with $df_p=0$ and
non-degenerate Hessian at $p$. Its tangent space at $p$
is $T_pX$, so the hypersurface Gauss equation gives
$(\GW)_p=R^X_p\neq0$, while its shape operator at $p$
is the Hessian of $f$, viewed as an endomorphism, for
the upward normal. After shrinking around $p$, $W$
remains non-degenerate. The same two ambient parallel
spinors still restrict to independent solutions.

Theorem~\ref{thm:A}, which is proved independently later,
shows that the number of common solutions in these
examples cannot exceed $N/2$. Thus each example has
exactly $N/2$ common solutions, proving
Theorem~\ref{thm:B}. The same count holds on sufficiently
small neighbourhoods of a point where the Gauss defect
is non-zero, so these also give local sharpness examples.
Ricci-flat metrics are real analytic in suitable
coordinates \cite[Theorem~5.2]{DeTurckKazdan1981};
the product, sphere and graph constructions can all
be taken real analytic.

\subsection{The Heisenberg group example}

Let $a>0$. Use the Heisenberg group law
\[
(u,v,w)*(x,y,z)=(u+x,v+y,w+z+auy)
\]
on $\R^3$, and set
\[
g = dx^2+dy^2+(dz-ax\,dy)^2,
\qquad
e_1=\partial_x,\quad
e_2=\partial_y+ax\,\partial_z,\quad
e_3=\partial_z.
\]
Left translation preserves the coframe $dx,dy,dz-ax\,dy$.
Consequently, both this metric and its dual orthonormal
frame $\{e_1,e_2,e_3\}$ are left-invariant. The frame
satisfies
\[
[e_1,e_2]=ae_3,\qquad
[e_1,e_3]=[e_2,e_3]=0.
\]
By the Koszul formula, $\nabla_{e_i}e_i=0$ for
$i=1,2,3$, and the non-zero connection coefficients are
\begin{align*}
\nabla_{e_1}e_2&=\frac a2e_3,
& \nabla_{e_1}e_3&=-\frac a2e_2,
& \nabla_{e_2}e_3&=\frac a2e_1,\\
\nabla_{e_2}e_1&=-\frac a2e_3,
& \nabla_{e_3}e_1&=-\frac a2e_2,
& \nabla_{e_3}e_2&=\frac a2e_1.
\end{align*}

Writing $\omega_{ij}(X)=g(\nabla_Xe_i,e_j)$, the spin
connection in the lifted frame is
\[
\nabla^{\Sigma M}_X\psi
=X(\psi)+\frac14\sum_{i,j}
\omega_{ij}(X)e_i\cdot e_j\cdot\psi,
\]
where $X(\psi)$ differentiates the spinor components.
We take the convention
$e_1\cdot e_2\cdot e_3=+\mathrm{id}$, which converts
the bivector actions below into vector actions.
For any spinor $\psi$ with constant components in this
frame, $X(\psi)=0$, and so
\begin{align}
\nabla^{\Sigma M}_{e_1}\psi
&= \frac a4e_2\cdot e_3\cdot\psi
= -\frac a4e_1\cdot\psi,
\nonumber\\
\nabla^{\Sigma M}_{e_2}\psi
&= -\frac a4e_1\cdot e_3\cdot\psi
= -\frac a4e_2\cdot\psi,
\label{eq:heisspin}\\
\nabla^{\Sigma M}_{e_3}\psi
&= -\frac a4e_1\cdot e_2\cdot\psi
= \frac a4e_3\cdot\psi.
\nonumber
\end{align}
Hence for
\begin{equation}
W = \frac a2\operatorname{diag}(1,1,-1)
\label{eq:heisW}
\end{equation}
with respect to $\{e_1,e_2,e_3\}$, every constant spinor
is a common solution. By injectivity of evaluation,
\[
\dim_\C\KW(M)=2=N.
\]
This is a concrete instance of the left-invariant
generalised Killing spinors treated in
\cite{Artacho2025}. With the opposite Clifford convention,
$W$ is replaced by $-W$. Writing the sectional curvature as
$K_{ij}=g(R(e_i,e_j)e_j,e_i)$, the connection coefficients
above give the first row of the following table.
The second row is the product of the corresponding
eigenvalues of $W$, as required by the Gauss equation.
\begin{tablelayout}
\[
\begin{tabular}{lccc}
\toprule
& $K_{12}$ & $K_{13}$ & $K_{23}$\\
\midrule
Actual curvature
& $-3a^2/4$ & $a^2/4$ & $a^2/4$\\
Required by the Euclidean Gauss equation
& $a^2/4$ & $-a^2/4$ & $-a^2/4$\\
\bottomrule
\end{tabular}
\]
\end{tablelayout}
Hence $\GW\neq0$, and Proposition~\ref{prop:C} follows.
A further direct computation gives
$\CW(e_1,e_2)=a^2e_3$. In this example the modified
connection itself is flat, but the Gauss and Codazzi
defects are both non-zero.

\section{Pointwise rigidity and the sharp smooth threshold}
\label{sec:smooth}

In this section we prove Theorem~\ref{thm:A}.
Fix $p\in M$ and use $E$, $\Eh$ and $\Vh_p$ from
Section~\ref{sec:conventions}. Throughout the proof
we assume
\[
\dim_\C\Vh_p=\dim_\C\KW(M)>N/2.
\]
By \eqref{eq:annihil}, every curvature bivector
$\Omega(X,Y)$ annihilates the same subspace $\Vh_p$
of the $N$-dimensional spin module.
We first bound the kernel of a non-zero bivector.
This shows that $\Omega(X,Y)=0$ except in dimensions
$7$, $11$ and $15$. In those cases, we classify the
bivectors that annihilate all of $\Vh_p$, then use the
curvature symmetries to show that $\Omega(X,Y)$ must
still vanish.

\subsection{Kernels of bivectors and parity-restricted
zero weights}

\subsubsection*{The eigenvalues of a bivector}
Set $q=\lfloor(n+1)/2\rfloor$, so the full spinor
module $\Delta_{n+1}$ of the formal ambient space
$\Eh$ has complex dimension $2^q$. A real bivector
$A\in\Lambda^2\Eh$ can be written, with respect to
a suitable orthonormal basis
$\{f_1,\dots,f_{n+1}\}$, as
\begin{equation}
A=\sum_{j=1}^q\mu_jf_{2j-1}\wedge f_{2j}.
\label{eq:normalform}
\end{equation}
When $n+1$ is even and the orientation of the basis
is fixed, we allow the $\mu_j$ to have either sign.
The Clifford operators $f_{2j-1}\cdot f_{2j}$ commute
with one another and square to $-\mathrm{id}$, so
they admit a basis of simultaneous eigenspinors.
For a non-zero simultaneous eigenspinor $\xi$, write
\[
(f_{2j-1}\cdot f_{2j})\cdot\xi
=i\varepsilon_j\xi,
\qquad \varepsilon_j\in\{\pm1\}.
\]
Every sign choice occurs. Indeed, multiplication by
$f_{2j-1}$ takes a simultaneous eigenspinor with signs
\[
(\varepsilon_1,\dots,\varepsilon_j,\dots,\varepsilon_q)
\quad\text{to one with signs}\quad
(\varepsilon_1,\dots,-\varepsilon_j,\dots,\varepsilon_q).
\]
This follows because $f_{2j-1}$ anticommutes with
$f_{2j-1}\cdot f_{2j}$ and commutes with the other
pairs. Since $f_{2j-1}^2=-1$, a non-zero spinor stays
non-zero. We can change any signs we wish in this
way, so there are $2^q$ simultaneous eigenspaces.
Their dimensions sum to $\dim_\C\Delta_{n+1}=2^q$,
and hence each has dimension one.

On the eigenspace with signs
$\varepsilon\in\{\pm1\}^q$, the action of $A$ is
multiplication by
\begin{equation}
i\sum_{j=1}^q\varepsilon_j\mu_j.
\label{eq:weights}
\end{equation}
It follows that $\dim_\C\ker(A\cdot)$ is exactly
the number of sign choices for which this sum
vanishes.

\begin{lemma}[Kernel bound on the full spin module]
\label{lem:fullkernel}
If $A\neq0$, then on the full spin module
$\Delta_{n+1}$,
\[
\dim_\C\ker(A\cdot)
\leq\frac12\dim_\C\Delta_{n+1}
=2^{q-1}.
\]
\end{lemma}

\begin{proof}
We may assume $\mu_1\neq0$. Once
$\varepsilon_2,\dots,\varepsilon_q$ are fixed, at most
one $\varepsilon_1\in\{\pm1\}$ satisfies the zero-weight
condition
$\varepsilon_1\mu_1=-\sum_{j=2}^q\varepsilon_j\mu_j$.
Hence there are at most $2^{q-1}$ zero weights.
\end{proof}

\subsubsection*{Even intrinsic dimension}
For even $n$, the formal ambient space has odd
dimension $n+1$. In this case, \eqref{eq:iota}
identifies the intrinsic spinor fibre with the full
ambient spin module $\Delta_{n+1}$, whose complex
dimension is $2^{n/2}=N$. For each fixed pair
$X,Y\in E$, \eqref{eq:annihil} shows that the kernel
of Clifford multiplication by $\Omega(X,Y)$ contains
$\Vh_p$, so its dimension exceeds $N/2$.
Lemma~\ref{lem:fullkernel} therefore forces
$\Omega(X,Y)=0$ for every $X,Y\in E$.

\subsubsection*{Odd intrinsic dimension}
Now suppose $n$ is odd, so $n=2q-1$. The formal
ambient dimension is $n+1=2q$, but \eqref{eq:iota}
identifies the intrinsic spinor fibre with just
one half-spin module:
\[
\dim_\C\Delta_{n+1}=2^q=2N,\qquad
\dim_\C\Delta^\tau_{n+1}=2^{q-1}=N,\qquad
N/2=2^{q-2}.
\]
Lemma~\ref{lem:fullkernel} bounds the kernel of a
non-zero bivector on the full module by $N$.
This does not contradict the existence of more than
$N/2$ independent solution values in the chosen half.
We therefore need a bound on the kernel within
that half-spin module.

On $\Delta^\tau_{n+1}$, the signs satisfy
$\prod_j\varepsilon_j=\sigma$,
where $\sigma\in\{\pm1\}$ depends on $\tau$
and our convention for the volume element. Setting
\[
Z^\sigma(\mu):=
\Bigl\{\varepsilon\in\{\pm1\}^q :
\prod_j\varepsilon_j=\sigma,\quad
\sum_j\varepsilon_j\mu_j=0\Bigr\},
\]
we have
$\dim_\C\ker(A\cdot|_{\Delta^\tau_{n+1}})
=|Z^\sigma(\mu)|$.

\subsubsection*{Solutions for $W$ and $-W$}
Before proving the general estimate, we illustrate
the weight calculation in dimension three.
Here $q=2$ and each ambient half-spin module has
complex dimension two. Suppose that there is a
non-zero generalised Killing spinor for each of
$W$ and $-W$. At a fixed point, the algebraic
version of \eqref{eq:iota-opp} is
$\iota^{\mathrm{opp}}(\xi)=e_0\cdot\iota(\xi)$.
Since
$\mathcal G_{-W}=\GW$ and
$\mathcal C_{-W}=-\CW$,
the opposite intertwining rule and
\eqref{eq:intid} give
\[
\iota^{\mathrm{opp}}
\bigl(R^{\nabla^{-W}}_{X,Y}\xi\bigr)
=\frac12\Omega(X,Y)\cdot
\iota^{\mathrm{opp}}(\xi).
\]
Thus the values of the two solutions give non-zero
vectors annihilated by the same curvature bivector
$\Omega(X,Y)$, one in each half-spin module.

For fixed $X,Y$, write
\[
\Omega(X,Y)
=\mu_1f_1\wedge f_2+\mu_2f_3\wedge f_4.
\]
The sign vectors $(+,+)$ and $(-,-)$ belong to one
parity class, and $(+,-)$ and $(-,+)$ to the other.
Thus, up to interchanging the half-spin labels,
the eigenvalues are
\[
\begin{aligned}
\text{on }\Delta_4^+ &: \quad
i(\mu_1+\mu_2),\; -i(\mu_1+\mu_2),\\
\text{on }\Delta_4^- &: \quad
i(\mu_1-\mu_2),\; -i(\mu_1-\mu_2).
\end{aligned}
\]
A non-zero kernel in each half therefore requires
$\mu_1+\mu_2=0$ and $\mu_1-\mu_2=0$, so both
coefficients vanish. Since $X,Y$ were arbitrary,
$\Omega=0$, and \eqref{eq:Omega-components} gives
the Gauss--Codazzi equations. This recovers the
Euclidean criterion of Lawn--Roth
\cite{LawnRoth2010}: one non-zero solution for each
sign suffices. The Heisenberg example in
Section~\ref{sec:examples} shows why both signs
matter: even two independent solutions for $W$
alone do not suffice.

In higher odd dimensions, requiring solutions for both
$W$ and $-W$ does not lower the threshold. Indeed, for
odd $n\geq5$, consider the example
\[
M=X^4\times S^{n-4}(r)
\subset X^4\times\R^{n-3}
\]
from Section~\ref{subsec:halfrank}. As computed there,
the ambient manifold has $N/2$ parallel spinors of
each chirality. Restriction through $\iota$ and
$\iota^{\mathrm{opp}}$ gives at least $N/2$ solutions
for $W$ and $-W$, respectively. Applying
Theorem~\ref{thm:A} to each field gives equality,
since $\mathcal G_{-W}=\GW$ and the Gauss equation
fails. Thus
\[
\dim_\C\KW(M)=\dim_\C\mathcal K_{-W}(M)=N/2.
\]
Even half the
possible solutions for each sign of $W$ do not
suffice to recover Euclidean hypersurface data.

\subsubsection*{The parity-restricted estimate}
We now return to the higher-dimensional problem,
where we must bound the kernel within a single
half-spin module.

\begin{theorem}[Parity-restricted zero-weight estimate]
\label{thm:LO}
For $q\geq2$, $\mu\in\R^q\setminus\{0\}$ and
$\sigma\in\{\pm1\}$,
\begin{equation}
|Z^\sigma(\mu)|\leq
\begin{cases}
2^{q-2}, & q \text{ odd},\\[2pt]
\max\Bigl\{2^{q-2},\dbinom{q}{q/2}\Bigr\},
& q\text{ even}.
\end{cases}
\label{eq:LObound}
\end{equation}
This bound is sharp for each $q$ and each fixed $\sigma$.
Moreover, if $q\in\{4,6,8\}$ and
$|Z^\sigma(\mu)|>2^{q-2}$, then all coefficients
are non-zero and
\begin{equation}
|\mu_1|=\dots=|\mu_q|,\qquad
|Z^\sigma(\mu)|=\binom{q}{q/2},\qquad
Z^{-\sigma}(\mu)=\varnothing.
\label{eq:LOrigid}
\end{equation}
\end{theorem}

The proof in Appendix~\ref{app:LO} uses the classical
antichain method \cite{Erdos1945,Sperner1928}.
We also prove the exceptional rigidity needed for the
geometric argument.

\subsubsection*{The exceptional dimensions}
If $\Omega(X,Y)\neq0$, its kernel on the chosen
half-spin module contains $\Vh_p$ and so has dimension
greater than $2^{q-2}$. For odd $q$,
Theorem~\ref{thm:LO} bounds that same dimension by
$2^{q-2}$, a contradiction. For even $q$, the only cases
in which the central binomial coefficient exceeds
$2^{q-2}$ are $q=2,4,6,8$. Since Theorem~\ref{thm:A}
assumes $n\geq4$, only $n=7,11,15$ remain.
The failure at $n=3$ was exhibited in
Proposition~\ref{prop:C}.

For the remaining dimensions $n=7,11,15$, we use the
metric to pass between two-forms, which are evaluated
on pairs of vectors, and bivectors, which act on
spinors by Clifford multiplication. For an orthogonal
complex structure $J$, its K\"ahler form and
corresponding bivector are characterised by
\[
\omega_J(X,Y):=\langle JX,Y\rangle
=\langle\omega_J^\sharp,X\wedge Y\rangle,
\qquad
\omega_J\in\Lambda^2\Eh^*,\quad
\omega_J^\sharp\in\Lambda^2\Eh.
\]
If all the coefficients in \eqref{eq:normalform} have
the same absolute value $\lambda>0$, then
\[
Jf_{2j-1}=\frac{\mu_j}{\lambda}f_{2j},\qquad
Jf_{2j}=-\frac{\mu_j}{\lambda}f_{2j-1}
\]
defines an orthogonal complex structure: each
$\mu_j/\lambda$ is $+1$ or $-1$, so
$J^2=-\mathrm{id}$ and $J$ preserves the metric.
Then $A=\lambda\omega_J^\sharp$. By the rigidity
statement in Theorem~\ref{thm:LO}, every non-zero
bivector whose kernel exceeds half rank in these
exceptional dimensions has this form.

\subsection{The common annihilator}

For $n=7,11,15$, the kernel bound alone is not enough.
Happily, there is another condition available: every
$\Omega(X,Y)$ must annihilate the same space
$\Vh_p$ of solution values. Accordingly,
let $q\in\{4,6,8\}$, so that
$\dim\Eh=2q$ and $N=2^{q-1}$, and retain the
hypothesis
$\dim_\C\Vh_p>N/2=2^{q-2}$.

\begin{definition}
The common annihilator of $\Vh_p$ is defined as
\[
\mathfrak a_p:=
\{A\in\Lambda^2\Eh :
A\cdot\xi=0\ \ \forall\xi\in\Vh_p\}.
\]
\end{definition}

This real vector space contains every curvature value
$\Omega(X,Y)$. If two bivectors annihilate $\Vh_p$,
so does their Clifford commutator.

\begin{lemma}
\label{lem:annihilator}
The common annihilator is one of
\[
\{0\},\qquad \R\omega_J^\sharp,\qquad
\operatorname{span}_\R
\{\omega_I^\sharp,\omega_J^\sharp,\omega_K^\sharp\}.
\]
In the last case, $I,J,K$ are orthogonal complex
structures forming a quaternionic triple with
$IJ=K=-JI$.
\end{lemma}

\begin{proof}
Bivectors have even degree and hence preserve the chosen
half-spin module. Using the commutator formula in
Lemma~\ref{lem:maximal} and
$[P,QR]=[P,Q]R+Q[P,R]$, we obtain
\begin{equation}
[f_if_j,f_kf_l]_{\Cl}
= 2\bigl(
\delta_{ik}f_j\wedge f_l
-\delta_{jk}f_i\wedge f_l
-\delta_{il}f_j\wedge f_k
+\delta_{jl}f_i\wedge f_k
\bigr).
\label{eq:comm2}
\end{equation}
Here we have substituted
$f_af_b=f_a\wedge f_b-\delta_{ab}$; the scalar terms
cancel. Hence $[A,B]_{\Cl}\in\Lambda^2\Eh$, and if
$A,B\in\mathfrak a_p$, then
\[
[A,B]_{\Cl}\cdot\xi
=A\cdot(B\cdot\xi)-B\cdot(A\cdot\xi)=0
\qquad(\xi\in\Vh_p).
\]
So $\mathfrak a_p$ is closed under the Clifford
commutator.

We now use the correspondence between bivectors and
skew endomorphisms from Lemma~\ref{lem:maximal} to
express the preceding rigidity as a condition on
their squares. To each bivector $A$ associate the
skew map $T_A$ defined by
\begin{equation}
\langle T_AX,Y\rangle
=\langle A,X\wedge Y\rangle.
\label{eq:TA}
\end{equation}
This correspondence is linear and injective, and
\[
T_{u\wedge v}z
=\langle u,z\rangle v-\langle v,z\rangle u,
\qquad
T_{\omega_J^\sharp}=J.
\]
Thus
\begin{equation}
[A,z]_{\Cl}=2T_Az,
\label{eq:Az}
\end{equation}
so the Lie algebra isomorphism
\eqref{eq:spin-adjoint} is precisely $A\mapsto2T_A$.
Since it preserves brackets,
\begin{equation}
T_{[A,B]_{\Cl}}=2[T_A,T_B].
\label{eq:Tbracket}
\end{equation}

If $0\neq A\in\mathfrak a_p$, its kernel contains
$\Vh_p$ and therefore exceeds half rank. By the rigidity
in Theorem~\ref{thm:LO} we can write
$A=\lambda\omega_J^\sharp$, so $T_A=\lambda J$ and
$T_A^2=-\lambda^2\mathrm{id}$. Here $\dim\Eh=2q$ and
$|A|^2=q\lambda^2$. Setting
\begin{equation}
Q(A):=-\frac1{2q}\tr(T_A^2)=\frac{|A|^2}{q},
\label{eq:Q}
\end{equation}
we have, for every $A\in\mathfrak a_p$ including zero,
\begin{equation}
T_A^2=-Q(A)\,\mathrm{id}.
\label{eq:TA2}
\end{equation}
Applying \eqref{eq:TA2} to $A+B$ and subtracting the
identities for $A$ and $B$ gives
\begin{equation}
T_AT_B+T_BT_A=-2b(A,B)\,\mathrm{id},
\qquad b(A,B)=\frac1q\langle A,B\rangle.
\label{eq:polar}
\end{equation}

If $\dim_\R\mathfrak a_p\geq2$, choose $A,B$
orthonormal with respect to $b$, and set $I=T_A$ and
$J=T_B$. Then $I^2=J^2=-\mathrm{id}$ and $IJ=-JI$, and
\[
K:=IJ=\frac12[I,J]=T_{\frac14[A,B]_{\Cl}}
\]
also corresponds to an element of $\mathfrak a_p$.
$K$ is skew and satisfies $K^2=-\mathrm{id}$, and
$I,J,K$ form a quaternionic triple. For example,
$IK+KI=0$, so by \eqref{eq:polar} the corresponding
bivectors are mutually $b$-orthogonal.

If there were a further non-zero element orthogonal
to this triple, its normalised associated map $L$
would anticommute with each of $I,J,K$.
But anticommutation with $I$ and $J$ also gives
$LK=LIJ=IJL=KL$. Hence $KL=0$, contradicting
$K^2=L^2=-\mathrm{id}$. The two-dimensional case does
not occur either, since two orthogonal generators
produce the third, so the possible dimensions are
$0$, $1$ and $3$. Applying \eqref{eq:TA2} to the
generator in the one-dimensional case as well gives
the stated classification.
\end{proof}

\subsection{Eliminating the exceptions by
curvature symmetries}

We continue to assume $q\in\{4,6,8\}$ and
$\dim_\C\Vh_p>2^{q-2}$. Set $h(X,Y)=g(W(X),Y)$.
By \eqref{eq:wedge} and the symmetry of $W$,
\[
g((W\wedge W)(X,Y)U,V)
=h(Y,U)h(X,V)-h(X,U)h(Y,V)
\]
has pair symmetry. Moreover, in the cyclic sum over
$X,Y,U$ all terms cancel by the symmetry of $h$.
Hence $\GW$, obtained by subtracting this term from $R$,
also satisfies
\begin{align}
\GW(X,Y,U,V)&=\GW(U,V,X,Y),
\label{eq:pair}\\
\GW(X,Y,U,V)+\GW(Y,U,X,V)+\GW(U,X,Y,V)&=0.
\label{eq:bianchi}
\end{align}

In what follows, all curvature arguments lie in
$E=e_0^\perp$. We use $\omega_J$ also for its restriction
to $E\times E$, while $\omega_J^\sharp$ continues to
denote the ambient bivector in $\Lambda^2\Eh$.
By \eqref{eq:Omega-components},
\[
\GW(X,Y,U,V)
=\langle\Omega(X,Y),U\wedge V\rangle.
\]

\begin{lemma}[The K\"ahler line case]
\label{lem:kahler}
If $\mathfrak a_p=\R\omega_J^\sharp$, then
$\Omega(X,Y)=0$ for all $X,Y\in E$.
\end{lemma}

\begin{proof}
Since $\Omega$ is alternating and bilinear in $X,Y$
and takes values in the fixed line
$\R\omega_J^\sharp$, its coefficient is a two-form
$\beta\in\Lambda^2E^*$. Thus
$\Omega(X,Y)=\beta(X,Y)\omega_J^\sharp$, and the
preceding identity gives
\[
\GW(X,Y,U,V)=\beta(X,Y)\omega_J(U,V).
\]
The subspace $H=\operatorname{span}\{e_0,Je_0\}$ and
its orthogonal complement are $J$-invariant. Choosing
a unit vector $u\in H^\perp$, the vectors $u$ and $Ju$
lie in $E$ and satisfy
$\omega_J(u,Ju)=\langle Ju,Ju\rangle=1$.
Setting $(U,V)=(u,Ju)$ in the pair symmetry gives,
for every $X,Y\in E$,
\[
\beta(X,Y)=\beta(u,Ju)\omega_J(X,Y).
\]
Hence $\beta=c\omega_J$ on $E$, where
$c=\beta(u,Ju)$, and
$\GW=c\,\omega_J\otimes\omega_J$.

The space $\operatorname{span}\{e_0,Je_0,u,Ju\}$
is $J$-invariant and has dimension four. Its
orthogonal complement is therefore $J$-invariant
and has dimension $2q-4\geq4$. Choose a unit vector
$v$ there. Then $u,Ju,v,Jv$ are orthonormal tangent
vectors. Applying \eqref{eq:bianchi} to this quadruple,
\[
0=c\bigl(
\omega_J(u,Ju)\omega_J(v,Jv)
+\omega_J(Ju,v)\omega_J(u,Jv)
+\omega_J(v,u)\omega_J(Ju,Jv)
\bigr)=c.
\]
The last two terms vanish by the orthogonality of
distinct $J$-invariant planes. Hence $c=0$, $\beta=0$
and $\Omega=0$.
\end{proof}

\begin{lemma}[The quaternionic triple case]
\label{lem:quaternionic}
If
$\mathfrak a_p=\operatorname{span}_\R
\{\omega_I^\sharp,\omega_J^\sharp,\omega_K^\sharp\}$,
then $\Omega(X,Y)=0$ for all $X,Y\in E$.
\end{lemma}

\begin{proof}
Set $(I_1,I_2,I_3)=(I,J,K)$ and
$\omega_a=\omega_{I_a}$, using the same restriction
convention as above. For some two-forms
$\beta_a\in\Lambda^2E^*$ we can write
\[
\Omega(X,Y)=\sum_a\beta_a(X,Y)\omega_a^\sharp,
\qquad
\GW(X,Y,U,V)=\sum_a\beta_a(X,Y)\omega_a(U,V).
\]
The quaternionic four-dimensional subspace
\[
H_0:=\operatorname{span}\{e_0,Ie_0,Je_0,Ke_0\}
\]
and its orthogonal complement are invariant under
each $I_a$. Here the dimension hypothesis provides
room for a second quaternionic four-plane: since
$\dim\Eh\geq8$, choose a unit vector $v\in H_0^\perp$.
Then $H_1=\operatorname{span}\{v,Iv,Jv,Kv\}$ is also
a quaternionic four-dimensional subspace, with
$H_1\subset E$ and $H_1\perp H_0$.
The quaternionic relations and orthogonality give
\begin{equation}
\omega_a(v,I_bv)=\langle I_av,I_bv\rangle=\delta_{ab}.
\label{eq:phiv}
\end{equation}
In particular, the restrictions of
$\omega_1,\omega_2,\omega_3$ to $E$ are linearly
independent.

Setting $(U,V)=(v,I_av)$ in the pair symmetry gives
\[
\beta_a(X,Y)=\sum_b\beta_b(v,I_av)\omega_b(X,Y).
\]
Defining $S_{ab}:=\beta_b(v,I_av)$, we get
$\beta_a=\sum_bS_{ab}\omega_b$. Evaluating this on
$(v,I_bv)$ and comparing with the definition gives
$S_{ab}=S_{ba}$. Hence
\begin{equation}
\GW(X,Y,U,V)
=\sum_{r,s}S_{rs}\,\omega_s(X,Y)\omega_r(U,V).
\label{eq:GS}
\end{equation}

If $x\in H_0\cap E$ and $y\in H_1$, then
$\omega_s(x,y)=\langle I_sx,y\rangle=0$, since
$I_sx\in H_0$. Equation~\eqref{eq:GS} therefore gives
$\GW(x,y,U,V)=0$ whenever $x\in H_0\cap E$,
$y\in H_1$ and $U,V\in E$. To show that every
$S_{ab}$ vanishes, choose $a,b\in\{1,2,3\}$
independently, and choose $c,d$ so that $(b,c,d)$
is one of $(1,2,3)$, $(2,3,1)$ or $(3,1,2)$.
By the quaternionic relations and \eqref{eq:phiv},
\[
\omega_s(I_ce_0,I_de_0)=\delta_{sb},\qquad
\omega_r(v,I_av)=\delta_{ra}.
\]
In \eqref{eq:GS}, the pair $(I_ce_0,I_de_0)$
selects column $b$, and the pair $(v,I_av)$ selects
row $a$. Only the term $S_{ab}$ remains.
All four vectors lie in $E$.
From \eqref{eq:bianchi} and \eqref{eq:GS},
\begin{align*}
0&=\GW(I_ce_0,I_de_0,v,I_av)\\
&\quad+\GW(I_de_0,v,I_ce_0,I_av)
+\GW(v,I_ce_0,I_de_0,I_av)\\
&=S_{ab}.
\end{align*}
The last two terms vanish because their first two
inputs lie in $H_0\cap E$ and $H_1$, respectively,
or in the reverse order. Since $a$ and $b$ can be
chosen independently, every $S_{ab}$ vanishes.
Hence $\beta_a=0$ and $\Omega=0$.
\end{proof}

\begin{proof}[Proof of Theorem~\ref{thm:A}]
At every $p\in M$,
$\dim_\C\Vh_p=\dim_\C\KW(M)>N/2$.
If $n$ is even, Lemma~\ref{lem:fullkernel} gives
$\Omega(X,Y)=0$. If $n$ is odd and
$n\notin\{3,7,11,15\}$, Theorem~\ref{thm:LO} gives
the same conclusion. For the remaining cases
$n=7,11,15$ we use Lemma~\ref{lem:annihilator}.
If the annihilator is zero, the conclusion is
immediate; if it is a line or a quaternionic triple,
we apply Lemmas~\ref{lem:kahler}
and~\ref{lem:quaternionic}. Hence $\Omega=0$ in all
cases. The component identities
\eqref{eq:Omega-components} give $\GW=\CW=0$. By the fundamental
theorem of hypersurfaces, the required local isometric
immersion exists.
\end{proof}

\begin{corollary}[Local multiplicity gap]
Let $(M^n,g)$ be a connected smooth Riemannian spin
manifold, $n\geq4$, and let $W$ be a smooth symmetric
endomorphism field of $TM$. For every connected,
simply connected open set $U\subset M$, the solution
space for $W|_U$ satisfies
\[
\dim_\C\KW(U)\in\{0,1,\dots,N/2\}\cup\{N\}.
\]
\end{corollary}

\begin{proof}
If $\dim_\C\KW(U)>N/2$, Theorem~\ref{thm:A} gives
$\GW=\CW=0$, so \eqref{eq:intid} gives
$R^{\nabla^W}=0$. Fix $p\in U$ and $\xi\in\Sigma_pM$.
For $x\in U$, define $\psi(x)$ by transporting $\xi$
from $p$ to $x$. Any two paths from $p$ to $x$ are
homotopic with fixed endpoints, since $U$ is simply
connected. Flatness then implies that they give the
same parallel transport. This defines the unique
$\nabla^W$-parallel section with $\psi(p)=\xi$.
Thus evaluation is an isomorphism
$\KW(U)\cong\Sigma_pM$, and $\dim_\C\KW(U)=N$.
\end{proof}

Without simple connectedness, the hypothesis
$\dim_\C\KW(U)>N/2$ still gives $R^{\nabla^W}=0$,
but parallel transport around a loop need not be the
identity. Fix $p\in U$. The monodromy representation
$\rho_p:\pi_1(U,p)\to\mathrm{GL}(\Sigma_pM)$ assigns
to each homotopy class the $\nabla^W$-parallel
transport around a representative loop based at $p$.
A spinor value at $p$ extends to a global parallel
section if and only if it returns to itself around
every loop. Evaluation therefore identifies
\[
\KW(U)\cong
\{\xi\in\Sigma_pM:\rho_p([\gamma])\xi=\xi
\text{ for every }[\gamma]\in\pi_1(U,p)\}.
\]
If the monodromy is non-trivial, this fixed space is
a proper subspace of $\Sigma_pM$. It follows that
$\dim_\C\KW(U)<N$. This is why simple connectedness
is needed for the conclusion above.

\section{Analytic extension and parallel spinor counts}
\label{sec:analytic}

Multiple common solutions constrain the holonomy of
an ambient extension only if they extend to the same
metric. For real analytic data, the metric extension
is determined by $(g,W)$, and the extension theorem
of Ammann--Moroianu--Moroianu \cite{AMM2013} applies
to every common solution on this same metric.
We can therefore use the holonomy classification to
bound the number of solutions when the ambient metric
is non-flat, keeping track of chirality in even
ambient dimension. This proves Theorem~\ref{thm:D}
and gives an independent proof of the sharp threshold
for real analytic data.

\subsection{The common Ricci-flat germ}

We first construct the ambient metric. Given a
connected Riemannian spin manifold $(M^n,g)$ and a
symmetric endomorphism field $W$, we want a Ricci-flat
metric $\bar g$ in one higher dimension such that
$M$ has induced metric $g$ and shape operator $W$.
Locally, take $Z=(-\varepsilon,\varepsilon)\times U$,
where $U\subset M$ is a sufficiently small contractible
neighbourhood. Identifying $U$ with the initial slice
$\{0\}\times U$, we use Gaussian coordinates:
\begin{equation}
\bar g=dt^2+g_t,\qquad
g_0=g|_U,\qquad
\dot g_0(X,Y)=-2g(W(X),Y).
\label{eq:cauchydata}
\end{equation}
Here $\partial_t$ is the unit normal to the slices,
and $t\mapsto(t,x)$ is a unit-speed normal geodesic.
The last identity ensures that the initial slice has
shape operator $W$, with our convention
$W(X)=-\bar\nabla_X\partial_t$.
If the ambient metric is Ricci-flat, contracting the
Gauss and Codazzi equations gives the constraints
\eqref{eq:constraints} below on $(g,W)$.
Write $\nabla^g$ for the Levi-Civita connection and
define the codifferential of $W$ by
\[
(\delta W)(X)
:=-\sum_i g((\nabla^g_{e_i}W)(e_i),X)
\]
for a local orthonormal frame $\{e_i\}$.
A single non-zero generalised Killing spinor supplies
these constraints; see the proof of
\cite[Corollary~3.1]{AMM2013}. We recall the calculation
using the Schr\"odinger--Lichnerowicz formula.

\begin{proposition}
If $\KW(M)\neq0$, then
\begin{equation}
d\tr(W)+\delta W=0,\qquad
\Scal_g+\tr(W^2)-(\tr(W))^2=0.
\label{eq:constraints}
\end{equation}
\end{proposition}

\begin{proof}
Take $0\neq\psi\in\KW(M)$. Since $M$ is connected,
$\psi$ has non-zero constant length. Write $\nabla$
for the spin connection in this proof. Symmetry of
$W$, Clifford contraction and the Leibniz rule give
\[
D\psi=\frac12\tr(W)\psi,
\qquad
D^2\psi=\frac12(d\tr(W))^\sharp\cdot\psi
+\frac14(\tr(W))^2\psi.
\]
Recall that the connection Laplacian is
\[
\nabla^*\nabla\psi
=-\sum_i\left(
\nabla_{e_i}\nabla_{e_i}\psi
-\nabla_{\nabla^g_{e_i}e_i}\psi
\right).
\]
Choosing a local orthonormal frame that is geodesic
at a point, and differentiating the generalised
Killing equation there, gives
\begin{align*}
\nabla^*\nabla\psi
&=\frac12\sum_i
  (\nabla^g_{e_i}W)(e_i)\cdot\psi
  -\frac14\sum_i W(e_i)\cdot W(e_i)\cdot\psi\\
&=-\frac12(\delta W)^\sharp\cdot\psi
  +\frac14\tr(W^2)\psi.
\end{align*}
Applying the Schr\"odinger--Lichnerowicz formula
$D^2=\nabla^*\nabla+\Scal_g/4$
\cite{Lichnerowicz1963} yields
\begin{equation}
\bigl(\Scal_g+\tr(W^2)-(\tr(W))^2\bigr)\psi
=2(d\tr(W)+\delta W)^\sharp\cdot\psi.
\label{eq:lich}
\end{equation}
Taking the real part of the Hermitian inner product
with $\psi$ shows that the scalar coefficient on the
left vanishes, since Clifford multiplication by real
vectors is skew-Hermitian and $|\psi|^2>0$.
Writing $V=(d\tr(W)+\delta W)^\sharp$, the remaining
equation is $V\cdot\psi=0$. Since
$V\cdot V=-|V|^2$, we obtain $V=0$.
This proves \eqref{eq:constraints}.
\end{proof}

Assume now that $(g,W)$ is real analytic and
$\KW(M)\neq0$. Koiso's Einstein extension theorem
\cite{Koiso1981}, in the form stated in
\cite[Theorem~2.1]{AMM2013} with Einstein constant zero,
provides a unique real analytic Ricci-flat metric germ
with initial data \eqref{eq:cauchydata}.
In Gaussian coordinates, the tangential Ricci-flat
equations determine $\ddot g_t$ in terms of $g_t$,
$\dot g_t$ and tangential derivatives.
With $g_0$ and $\dot g_0$ as initial data, the
Cauchy--Kowalevski theorem gives a unique solution
$g_t$ for sufficiently small $|t|$. The constraints
\eqref{eq:constraints} are preserved along the
solution, so
\[
\operatorname{Ric}_{\bar g}(\partial_t,\cdot)=0.
\]
Together with the tangential equations, this gives
$\operatorname{Ric}_{\bar g}=0$. In fixed Gaussian coordinates, two metrics define the same \emph{germ} if they agree on a smaller
neighbourhood of the initial slice. Intrinsically, any two real analytic Ricci-flat
extensions $(Z_1,\bar g_1)$ and $(Z_2,\bar g_2)$
with the same initial data are related, after
shrinking around $U$, by a diffeomorphism
$\Phi:Z_1\to Z_2$ satisfying
\[
\Phi^*\bar g_2=\bar g_1,\qquad
\Phi|_U=\mathrm{id}_U,\qquad
d\Phi(\nu_1)=\nu_2,
\]
where $\nu_i$ is the chosen unit normal along the
initial slice in $Z_i$. We may thus work with
a metric $\bar g$ on
$Z=(-\varepsilon,\varepsilon)\times U$, shrinking
$U$ and $\varepsilon$ as needed. Its construction
uses only $g$ and $W$, so we keep $\bar g$ fixed
as $\psi$ varies in $\KW(U)$.

The spin structure on $U$, together with the
normal vector $\partial_t$, gives the product
spin structure on $Z$. Along the initial slice,
the identifications used in \eqref{eq:spingauss}
are
\[
\Sigma U\cong
\begin{cases}
(\Sigma Z)|_U,&n\text{ even},\\
(\Sigma^\tau Z)|_U,&n\text{ odd}.
\end{cases}
\]
For $\psi\in\KW(U)$, the spinorial Gauss formula
then gives
\[
\bar\nabla^{\Sigma Z}_X\psi
=\nabla^{\Sigma U}_X\psi
+\frac12W(X)\cdot\psi
=0,\qquad X\in TU.
\]
Parallel transport along the normal geodesics gives
the unique candidate for an ambient parallel extension:
\[
\bar\nabla^{\Sigma Z}_{\partial_t}\Psi=0,\qquad
\Psi|_{\{0\}\times U}=\psi.
\]
The parallel spinor extension theorem of
Ammann--Moroianu--Moroianu
\cite[Theorem~3.2]{AMM2013} shows that this extension
is parallel in every direction. Here the initial
slice is oriented and real analytic, and
$(Z,\bar g)$ is Ricci-flat and spin, as required.
The theorem also requires
$i_*:\pi_1(U)\to\pi_1(Z)$ to be surjective,
where $i:U\hookrightarrow Z$ is the inclusion.
In our product neighbourhood $i_*$ is an
isomorphism: $Z$ deformation retracts onto $U$.
Every $\psi\in\KW(U)$ therefore extends to a
parallel spinor on the same $Z$. The extension is
linear, as it is defined by parallel transport.
For odd $n$, the ambient spin connection preserves
$\Sigma^\tau Z$, so the extension lies in this
chosen half-spin bundle.
By the spinorial Gauss formula~\eqref{eq:spingauss},
the restriction of any parallel spinor in the chosen
ambient bundle belongs to $\KW(U)$.
Restriction therefore gives
\begin{equation}
\begin{aligned}
\mathcal P(Z)&\xrightarrow{\ \cong\ }\KW(U),
&&n\text{ even},\\
\mathcal P^\tau(Z)&\xrightarrow{\ \cong\ }\KW(U),
&&n\text{ odd}.
\end{aligned}
\label{eq:restrictioniso}
\end{equation}
Surjectivity follows from the extension theorem, and
injectivity was established in
Section~\ref{sec:examples}.

\subsection{The non-flat upper bound and its attainment}

We now count the parallel spinors on the common
Ricci-flat extension $Z^{n+1}$ to bound the number of
solutions for the fixed field $W$. Here $n$ is the
hypersurface dimension, $N=2^{\lfloor n/2\rfloor}$ is
the intrinsic spinor rank, and $d$ denotes the dimension
of an individual factor.
On a sufficiently small simply connected domain,
evaluation at $p\in Z$ identifies parallel spinors
with the vectors fixed by the lifted holonomy:
\[
\mathcal P(Z)
\xrightarrow[\Psi\mapsto\Psi(p)]{\ \cong\ }
\left\{\xi\in\Sigma_pZ:
h\cdot\xi=\xi\
\text{for all }h\in\widetilde{\Hol}_p(Z)\right\}.
\]
Here $\widetilde{\Hol}_p(Z)$ is the holonomy of the
spin connection. A fixed vector extends by parallel
transport independently of the chosen path, and
every parallel spinor is determined by its value at
$p$. When $n$ is odd, $Z^{n+1}$ has even dimension,
and the same correspondence holds in each half-spin
bundle.

Since $\bar g$ is real analytic, its restricted
holonomy is unchanged on connected neighbourhoods
of $p$
\cite[Chapter~II, Theorem~10.8]{KobayashiNomizu1963}.
After shrinking around $p$, the local de Rham
decomposition gives
\[
Z\cong\R^\ell\times Z_1\times\cdots\times Z_k,
\]
where the $Z_i$ are non-flat irreducible factors.
Each factor is Ricci-flat, since the Ricci tensor
splits under products, and each has parallel spinors
by the product rule below. A Ricci-flat locally
symmetric Riemannian space is flat, so the non-flat
factors are not locally symmetric. The classification
of irreducible restricted holonomy applies locally
to these factors. Here $\Hol^0$ denotes the identity
component of the Riemannian holonomy group. On our
simply connected domains, this is the full holonomy
group. Parallel spinors correspond to vectors fixed
by its lift to the spin group.
For each factor $Y$, write $d=\dim Y$ and
$p(Y)=\dim_\C\mathcal P(Y)$ for its number of parallel
spinors. In even dimensions, write
$p^\pm(Y)=\dim_\C\mathcal P^\pm(Y)$ for the number of
each chirality; reversing orientation may interchange
these counts. The table below gives Wang's counts
\cite{Wang1989}, with the holonomy lifts discussed in
\cite{MoroianuSemmelmann2000}. We take $r\geq2$ for
$\mathrm{SU}(r)$ and $r\geq1$ for $\mathrm{Sp}(r)$;
the case $\mathrm{SU}(2)=\mathrm{Sp}(1)$ appears
in both families.
\begin{tablelayout}
\setlength{\arraycolsep}{10pt}
\begin{equation}
\begin{array}{lcccc}
\toprule
\Hol^0 & d & p^+ & p^- & p\\
\midrule
\mathrm{SU}(r),\ r\text{ even}
& 2r & 2 & 0 & 2\\
\mathrm{SU}(r),\ r\text{ odd}
& 2r & 1 & 1 & 2\\
\mathrm{Sp}(r) & 4r & r+1 & 0 & r+1\\
G_2 & 7 & \text{--} & \text{--} & 1\\
\mathrm{Spin}(7) & 8 & 1 & 0 & 1\\
\bottomrule
\end{array}
\label{eq:wangtable}
\end{equation}
\end{tablelayout}

Following Alekseevsky--Cort\'es
\cite[Section~3]{AlekseevskyCortes2009}, we consider
the proportion of the full spinor fibre represented
by parallel spinors. For a factor $Y$ of dimension
$d$, define its \emph{full-spin efficiency} by
\begin{equation}
\eta(Y):=\frac{p(Y)}{2^{\lfloor d/2\rfloor}}.
\label{eq:eta}
\end{equation}
From \eqref{eq:wangtable}, the efficiencies of
the non-flat irreducible factors are
\[
\eta_{\mathrm{SU}(r)}=\frac2{2^r},\qquad
\eta_{\mathrm{Sp}(r)}=\frac{r+1}{2^{2r}},\qquad
\eta_{G_2}=\frac18,\qquad
\eta_{\mathrm{Spin}(7)}=\frac1{16}.
\]
All are at most $1/2$, with equality precisely in
the four-dimensional $\mathrm{Sp}(1)$ case.
A flat factor has efficiency $1$.

Alekseevsky--Cort\'es prove that this invariant is
multiplicative under products
\cite[Corollary~3]{AlekseevskyCortes2009}.
We recall the product formulas underlying their
argument \cite[Corollary~2]{AlekseevskyCortes2009}.
For factors $Y_1,Y_2$ of dimensions $d_1,d_2$,
the full spin representation restricts to the product
of the spin groups as
\[
\Delta_{d_1+d_2}\cong
\begin{cases}
\Delta_{d_1}\otimes\Delta_{d_2},
  &\text{at least one }d_i\text{ is even},\\
(\Delta_{d_1}\otimes\Delta_{d_2})^{\oplus2},
  &\text{both }d_i\text{ are odd}.
\end{cases}
\]
The two holonomy groups act on their respective tensor
factors. Thus, writing $G_i$ for the lifted holonomy
of $Y_i$,
\[
(\Delta_{d_1}\otimes\Delta_{d_2})^{G_1\times G_2}
=\Delta_{d_1}^{G_1}\otimes\Delta_{d_2}^{G_2}.
\]
The corresponding rule for parallel spinors is
\[
\mathcal P(Y_1\times Y_2)\cong
\begin{cases}
\mathcal P(Y_1)\otimes\mathcal P(Y_2),
  &\text{at least one }d_i\text{ is even},\\
(\mathcal P(Y_1)\otimes\mathcal P(Y_2))^{\oplus2},
  &\text{both }d_i\text{ are odd}.
\end{cases}
\]
The number of parallel spinors and the rank of the
full spinor bundle therefore acquire the same factor
$c$, where $c=2$ if both dimensions are odd and $c=1$
otherwise. Consequently,
\[
\eta(Y_1\times Y_2)
=\frac{c\,p(Y_1)p(Y_2)}
 {c\,2^{\lfloor d_1/2\rfloor}
       2^{\lfloor d_2/2\rfloor}}
=\eta(Y_1)\eta(Y_2).
\]
For a non-flat germ, the de Rham decomposition has
$k\geq1$, so
$\eta(Z)=\prod_{i=1}^k\eta(Z_i)\leq1/2$.
This gives the full-spin bound
\begin{equation}
p(Z)\leq\frac12\,2^{\lfloor(n+1)/2\rfloor}.
\label{eq:fullbound}
\end{equation}

The use of Wang's counts and the de Rham decomposition
to obtain a full-spin flatness criterion already
appears in Alekseevsky--Cort\'es
\cite[proof of Theorem~4(ii)]{AlekseevskyCortes2009}.
Their discussion of equality also identifies the
product of a flat manifold and a manifold with
holonomy $\mathrm{SU}(2)$, which we use below to
attain the bounds.

When $n$ is even, $Z^{n+1}$ has odd dimension and its
full spinor bundle has rank $2^{n/2}=N$, equal to
the intrinsic spinor rank. Thus
\eqref{eq:restrictioniso} already gives
$\dim_\C\KW(U)=p(Z)\leq N/2$.
When $n$ is odd, the correspondence instead uses one
half-spin space $\mathcal P^\tau(Z)$, and
\[
\operatorname{rank}_\C\Sigma Z=2^{(n+1)/2}=2N,\qquad
\operatorname{rank}_\C\Sigma^\pm Z=2^{(n-1)/2}=N.
\]
The bound \eqref{eq:fullbound} only gives
$p^+(Z)+p^-(Z)\leq N$, whereas we need
$p^\tau(Z)\leq N/2$. We therefore need a bound for
each chirality separately.
For two even-dimensional factors $Y_1,Y_2$, the
combinations $(+,+)$ and $(-,-)$ give positive
chirality on the product; $(+,-)$ and $(-,+)$ give
negative chirality.
Writing $p_i^\pm=p^\pm(Y_i)$, the product rules are
\begin{equation}
\begin{aligned}
p^+(Y_1\times Y_2)&=p_1^+p_2^++p_1^-p_2^-,\\
p^-(Y_1\times Y_2)&=p_1^+p_2^-+p_1^-p_2^+.
\end{aligned}
\label{eq:chiralproduct}
\end{equation}
In particular, if $p_1^+=p_1^-=a$, both product
counts equal $a(p_2^++p_2^-)$. Thus, if either factor
has equally many parallel spinors of each chirality,
so does the product.
If both factors are odd-dimensional, multiplication
by either factor's volume form takes positive
parallel spinors on the product to negative ones,
and conversely. It preserves parallelness because
the volume form is parallel. Its Clifford action
is invertible, so the two spaces have equal dimension.
The full-spin product rule gives a total of
$2p(Y_1)p(Y_2)$ parallel spinors, so
\[
p^+(Y_1\times Y_2)=p^-(Y_1\times Y_2)=p(Y_1)p(Y_2).
\]
We want to bound each chiral count by one quarter of
the full spinor rank. If the two counts are equal,
the preceding bound on their sum immediately gives
this. Otherwise, we will show that the total number
of parallel spinors is itself at most one quarter
of the full spinor rank.

\begin{proposition}[Half-spin upper bound]
\label{prop:halfbound}
Let $n\geq5$ be odd. For a non-flat Ricci-flat germ
$Z^{n+1}$ with non-zero parallel spinors,
\begin{equation}
\max\{p^+(Z),p^-(Z)\}\leq2^{(n-3)/2}=N/2.
\label{eq:halfbound}
\end{equation}
\end{proposition}

\begin{proof}
Since $\dim Z=n+1$ is even, the local product
decomposition has an even number of odd-dimensional
factors, including the flat factor when its dimension
is odd. Pair these factors and keep the
even-dimensional factors separately.
Suppose first that one factor or pair has equal
numbers of parallel spinors of each chirality.
By the product rules, the same is true of $Z$, and
\eqref{eq:fullbound} gives
\[
p^\pm(Z)=\frac12p(Z)
\leq\frac14\,2^{(n+1)/2}
=2^{(n-3)/2}=N/2.
\]
This case includes any flat factor of positive even
dimension, any $\mathrm{SU}(r)$ factor with $r$ odd,
and any pair of odd-dimensional factors.

Suppose now that no factor or pair has equal chiral
counts. By \eqref{eq:wangtable}, all the factors are
then even-dimensional and non-flat, with parallel
spinors of only one chirality. If there are at least
two factors, multiplicativity gives
$\eta(Z)\leq(1/2)^2=1/4$. If there is only one factor,
the four-dimensional $\mathrm{Sp}(1)$ case is
excluded since $\dim Z=n+1\geq6$. Every remaining
entry of \eqref{eq:wangtable} with parallel spinors
of only one chirality has full-spin efficiency at
most $1/4$. In either case,
\[
p^\pm(Z)\leq p(Z)=\eta(Z)\,2^{(n+1)/2}
\leq2^{(n-3)/2}=N/2.\qedhere
\]
\end{proof}

These bounds are attained by taking
$Z^{n+1}=X^4\times\R^{n-3}$, where $X$ is the non-flat
hyperk\"ahler factor of Section~\ref{subsec:halfrank},
oriented so that $p^+(X)=2$ and $p^-(X)=0$.
For even $n\geq4$, the equality
$\eta(Z)=\eta(X)=1/2$ gives $p(Z)=N/2$, attaining
\eqref{eq:fullbound}.
For odd $n\geq5$, the flat factor has positive even
dimension, with $p^\pm(\R^{n-3})=2^{(n-5)/2}$, so
\eqref{eq:chiralproduct} gives
\begin{equation}
p^+(Z)=p^-(Z)
=2\cdot2^{(n-5)/2}=2^{(n-3)/2}=N/2,
\label{eq:attained}
\end{equation}
attaining \eqref{eq:halfbound}. For $n=7$, there is
also an equality example without a flat factor:
if $X_1^4,X_2^4$ are non-flat hyperk\"ahler factors
oriented as above, then
$p^+(X_1\times X_2)=4$ and
$p^-(X_1\times X_2)=0$.

For $n=3$, however, we may take $Z=X^4$, with
$p^+(Z)=2=\operatorname{rank}_\C\Sigma^+Z=N$.
The entire positive half-spin fibre is realised by
parallel spinors despite the metric being non-flat.
This explains why maximal generalised Killing spinor
multiplicity on a three-dimensional hypersurface
does not force the ambient metric to be flat.

\begin{proof}[Completion of the proof of
Theorem~\ref{thm:D}]
The simultaneous extension and the correspondence
were proved in \eqref{eq:restrictioniso}.
For a non-flat germ, if $n$ is even,
\eqref{eq:fullbound} gives
$\dim_\C\KW(U)=p(Z)\leq2^{n/2-1}=N/2$.
If $n$ is odd and $n\geq5$,
Proposition~\ref{prop:halfbound} gives
$\dim_\C\KW(U)=p^\tau(Z)\leq2^{(n-3)/2}=N/2$.
The equality examples above and the hypersurfaces of
Section~\ref{subsec:halfrank} show that the bound is
attained for every $n\geq4$.
\end{proof}

\begin{corollary}[Forcing for analytic data]
Let $n\geq4$ and suppose that $(g,W)$ is real analytic.
If $\dim_\C\KW(U)>N/2$ on a sufficiently small connected,
simply connected domain $U\subset M$, then the common
Ricci-flat germ is flat. Hence $(g,W)$ satisfies the
Euclidean Gauss--Codazzi equations on $U$.
\end{corollary}

\begin{proof}
If the germ were non-flat, its number of parallel
spinors would contradict the upper bound in
Theorem~\ref{thm:D}. Applying the Gauss and Codazzi
formulae to the hypersurface in the flat germ gives
the conclusion.
\end{proof}

\appendix

\section{The parity-restricted zero-weight estimate}
\label{app:LO}

Here we prove Theorem~\ref{thm:LO}. The problem is
to estimate, for $\mu\in\R^q\setminus\{0\}$, the
number of $\varepsilon\in\{\pm1\}^q$ satisfying
simultaneously
\[
\sum_{j=1}^q\varepsilon_j\mu_j=0,\qquad
\prod_{j=1}^q\varepsilon_j=\sigma,
\]
and to determine the coefficients when this number
is exceptionally large.

\subsection{The antichain estimate}

A family $\mathcal A$ of subsets of
$[q]=\{1,\dots,q\}$ is an antichain if no two distinct
members are related by inclusion.
The Lubell--Yamamoto--Meshalkin inequality
\cite{Lubell1966,Yamamoto1954,Meshalkin1963} gives
\begin{equation}
\sum_{A\in\mathcal A}\binom{q}{|A|}^{-1}\leq1.
\label{eq:LYM}
\end{equation}
For the proof we refer to the works cited. Since
each binomial coefficient is at most the central
binomial coefficient, Sperner's bound
\cite{Sperner1928}
\begin{equation}
|\mathcal A|\leq\binom{q}{\lfloor q/2\rfloor}
\label{eq:sperner}
\end{equation}
follows immediately. Indeed, the left-hand side of
\eqref{eq:LYM} is at least
$|\mathcal A|/\binom{q}{\lfloor q/2\rfloor}$.

\subsection{Changing signs and zero coefficients}

For each $j$ write $\mu_j=s_ja_j$, with
$a_j=|\mu_j|$ and $s_j\in\{\pm1\}$; if $\mu_j=0$,
set $s_j=1$. Setting $\varepsilon'_j=s_j\varepsilon_j$,
\begin{equation}
\sum_j\varepsilon_j\mu_j=\sum_j\varepsilon'_ja_j,
\qquad
\prod_j\varepsilon'_j=\sigma\prod_js_j.
\label{eq:signchange}
\end{equation}
This is a bijection of sign vectors, so the signs
of the coefficients can be removed at the cost,
if necessary, of interchanging the two parity classes.

We first treat the case in which some coefficient
vanishes. After relabelling we may assume
$\mu_1\neq0$ and $\mu_q=0$. Once
$\varepsilon_2,\dots,\varepsilon_{q-1}$ are fixed,
the zero-sum condition determines $\varepsilon_1$
in at most one way. Once it is determined, the parity
condition determines $\varepsilon_q$ uniquely.
Hence
\begin{equation}
|Z^\sigma(\mu)|\leq2^{q-2}
\qquad\text{(when some coefficient is zero).}
\label{eq:zerocoeff}
\end{equation}
This case has to be treated separately, because the
antichain property of positive weights used later
may fail when there are zero coefficients.

\subsection{The sharp zero-weight upper bound}

Suppose that all coefficients are non-zero.
Set $a_j=|\mu_j|>0$ and
$\sigma'=\sigma\prod_j s_j$. The sign change
\eqref{eq:signchange} gives
$|Z^\sigma(\mu)|=|Z^{\sigma'}(a)|$.
Set $T=\sum_ja_j$. Associating to a sign vector
$\varepsilon'$ its set of positive indices
$A=\{j:\varepsilon'_j=1\}$, we have
\begin{equation}
\sum_j\varepsilon'_ja_j=0
\iff\sum_{j\in A}a_j=T/2,\qquad
\prod_j\varepsilon'_j=(-1)^{q-|A|}.
\label{eq:subsets}
\end{equation}
Hence prescribing the parity class prescribes the
parity of $|A|$. Let $\mathcal F$ be the family of
all subsets of weight $T/2$. If $A\subsetneq B$,
positivity gives
$\sum_{j\in A}a_j<\sum_{j\in B}a_j$,
so $\mathcal F$ is an antichain. Moreover, if
$A\in\mathcal F$ then its complement $A^c$ also
belongs to $\mathcal F$.

When $q$ is odd, complementation interchanges the
parity of the cardinality. Hence the numbers of
zero weights in the two parity classes are equal,
and \eqref{eq:sperner} gives
\[
|Z^{\sigma'}(a)|=\frac12|\mathcal F|
\leq\frac12\binom{q}{(q-1)/2}\leq2^{q-2}.
\]
The last inequality holds because, for odd $q$,
there are two equal largest binomial coefficients,
whose sum is at most $2^q$.

When $q$ is even, the subfamily
$\mathcal A\subset\mathcal F$ of the chosen parity
is itself an antichain, so
\[
|Z^{\sigma'}(a)|
=|\mathcal A|\leq\binom{q}{q/2}.
\]
Combining this with \eqref{eq:zerocoeff} gives
\eqref{eq:LObound}.

We check sharpness. For $q\geq3$ take
$\mu=(1,1,0,\dots,0)$. The zero-sum condition is
$\varepsilon_1=-\varepsilon_2$.
Since the product
of the remaining $q-2$ signs is prescribed, each
parity class contains $2\cdot2^{q-3}=2^{q-2}$
solutions. On the other hand, for even $q$ take
$\mu=(1,\dots,1)$. The zero weights are then exactly
those with $q/2$ positive and $q/2$ negative signs,
and there are $\binom{q}{q/2}$ of them. They all
belong to the same parity class. Changing the sign
of one coefficient interchanges the two parity
classes. Sharpness for $q=2$ is also given by the
latter example. Hence the upper bound is attained
for each fixed $\sigma$.

We determine for which even $q$ the central binomial
coefficient exceeds half rank.
\begin{tablelayout}
\[
\begin{tabular}{cccccc}
\toprule
$q$ & 2 & 4 & 6 & 8 & 10\\
\midrule
$\binom{q}{q/2}$ & 2 & 6 & 20 & 70 & 252\\
$2^{q-2}$ & 1 & 4 & 16 & 64 & 256\\
\bottomrule
\end{tabular}
\]
\end{tablelayout}
Setting $c_q=\binom{q}{q/2}/2^{q-2}$, we have
\[
\frac{c_{q+2}}{c_q}=\frac{q+1}{q+2}<1,
\]
so the exceptions are exactly $q=2,4,6,8$.

\subsection{Rigidity in the exceptional cases}

We first prove a lemma forcing the coefficients
to have equal size.

\begin{lemma}[Middle-layer swap lemma]
\label{lem:swap}
Let $q=2m$, $a_j>0$ and $T=\sum_ja_j$. Suppose that
each member of a family $\mathcal A$ of $m$-element
subsets has weight $T/2$. If fewer than
$\binom{q-2}{m-1}$ of all $m$-element subsets are
missing from $\mathcal A$, then $a_1=\dots=a_q$.
\end{lemma}

\begin{proof}
Fix distinct $i,j$, and for
$B\subset[q]\setminus\{i,j\}$ with $|B|=m-1$, pair
\[
\{i\}\cup B,\qquad\{j\}\cup B.
\]
There are $\binom{q-2}{m-1}$ such pairs, and they
are pairwise disjoint. Fewer than
$\binom{q-2}{m-1}$ subsets are missing, so at least
one pair has both members in $\mathcal A$.
Both have weight $T/2$ and share the same subset $B$.
Cancelling the weights of the elements of $B$ gives
$a_i=a_j$. This holds for every pair $i\ne j$, so
$a_1=\dots=a_q$.
\end{proof}

Assume $q\in\{4,6,8\}$ and
$|Z^\sigma(\mu)|>2^{q-2}$. By \eqref{eq:zerocoeff}
all coefficients are non-zero, and after changing
signs we may assume they are all positive.
Let $\mathcal A$ be the family of zero-weight
subsets of the prescribed parity, and set
\[
b_k:=|\{A\in\mathcal A:|A|=k\}|.
\]
Since $q$ is even, complementation preserves this
family, and $b_k=b_{q-k}$. Moreover, complementation
has no fixed points, so $|\mathcal A|$ is even.
The empty set and the whole set do not have weight
$T/2$. If a singleton $\{i\}$ belongs to the family,
so does its complement. Any other member either
contains $i$, and so properly contains $\{i\}$,
or does not contain $i$, and so is properly contained
in the complement. By the antichain property, in
this case $|\mathcal A|=2$.

\subsubsection*{The case $q=4$}
If the chosen parity is odd, the possible sizes are
$1$ and $3$, and the LYM inequality gives
\[
\tfrac{b_1}{4}+\tfrac{b_3}{4}\leq1,\qquad
|\mathcal A|\leq4,
\]
contrary to the hypothesis. Hence every member has
size $2$. There are $6$ such subsets in total, and
since $|\mathcal A|$ is even and greater than $4$,
we get $|\mathcal A|=6$. No middle-layer subsets
are missing, and since $0<\binom21=2$,
Lemma~\ref{lem:swap} applies.

\subsubsection*{The case $q=6$}
If the chosen parity is even, the possible sizes are
$2$ and $4$, and
\[
\tfrac{b_2}{15}+\tfrac{b_4}{15}\leq1
\]
gives $|\mathcal A|\leq15<16$. Hence the parity is
odd, and the sizes are $1$, $3$ and $5$.
The observation on singletons excludes sizes $1$
and $5$, so every member has size $3$.
Since $\binom63=20$, $|\mathcal A|>16$, and
$|\mathcal A|$ is even, at most two subsets are
missing. Since $2<\binom42=6$,
Lemma~\ref{lem:swap} applies.

\subsubsection*{The case $q=8$}
If the odd parity is chosen and a singleton is a
member, there are two members. Otherwise only
the sizes $3$ and $5$ are possible, and
\[
\tfrac{b_3}{56}+\tfrac{b_5}{56}\leq1
\]
gives $|\mathcal A|\leq56<64$. Hence the parity is
even, and the sizes are $2$, $4$ and $6$.
From $b_2=b_6$ and the LYM inequality,
\[
\tfrac{2b_2}{28}+\tfrac{b_4}{70}\leq1,\qquad
b_4\leq70-5b_2,\qquad
|\mathcal A|\leq70-3b_2.
\]
If $b_2\geq2$, the count is at most $64$.
Suppose $b_2=1$, and let $P$ be the unique
two-element member of $\mathcal A$. Its complement
$P^c$ is the unique six-element member.
A four-element member $A$ cannot contain $P$,
since $\mathcal A$ is an antichain. Nor can $A$
be disjoint from $P$: then $A\subset P^c$, again
contradicting the antichain property.
Thus $|A\cap P|=1$. There are at most
$\binom21\binom63=40$ such subsets, so
\[
|\mathcal A|\leq1+40+1=42<64.
\]
Hence $b_2=b_6=0$, and every member has size $4$.
Since $\binom84=70$ and $|\mathcal A|$ is even,
$|\mathcal A|\geq66$, so at most four subsets are
missing. Since $4<\binom63=20$,
Lemma~\ref{lem:swap} applies.

\medskip

In all three cases we have obtained
$a_1=\dots=a_q$. Restoring the original signs gives
$|\mu_1|=\dots=|\mu_q|\neq0$. After the sign change,
the zero-sum condition is equivalent to there being
exactly $q/2$ positive signs. Hence there are
$\binom{q}{q/2}$ zero weights, all in the same parity
class, and none in the opposite class.
This completes the proof of Theorem~\ref{thm:LO}.

\section*{AI use}

Large language models were used during the
research and preparation of this article to explore
proof strategies, derive and check formulas, examine
arguments for gaps, and assist with exposition.
In particular, they contributed substantially to developing the arguments for the exceptional dimensions $7$, $11$
and $15$ in the proof of the smooth threshold theorem.
The authors take full responsibility for the content
and correctness of the article.

\bibliographystyle{amsalpha}
\bibliography{references}

\end{document}